\documentclass[a4paper, 11pt, twoside, dvips]{article}
\usepackage{amsmath}
\usepackage{amsfonts}
\usepackage{float}
\usepackage{amssymb}
\usepackage{amsthm}
\usepackage{amsfonts}
\usepackage{verbatim}
\usepackage{graphicx}
\usepackage{psfrag}
\usepackage[T1]{fontenc}
\usepackage{fancyhdr}
\theoremstyle{plain}
\newtheorem{thm}{Theorem}[section]
\newtheorem{prop}[thm]{Proposition}
\newtheorem{lemma}[thm]{Lemma}
\newtheorem{defn}[thm]{Definition}
\newtheorem{cor}[thm]{Corollary}

\newtheorem{rem}[thm]{Remark}

\newcommand{\df}{\stackrel{\mathrm{def}}{=}}
\newcommand{\R}[1]{\mathbb{R}^{#1}}
\newcommand{\mE}{\mathbb{E}}
\newcommand{\mP}{\mathbb{P}}
\newcommand{\om}{\omega}
\newcommand{\all}{\om \in \Omega}
\newcommand{\cL}{\mathcal L}
\newcommand{\cF}{\mathcal F}
\newcommand{\cN}{\mathcal N}

\numberwithin{equation}{section}

\begin{document}
\title{Examples of Condition $(T)$ for Diffusions in a Random Environment }
\author{Tom Schmitz\\
        Department of Mathematics\\
        ETH Zurich\\
        CH-8092 Zurich\\
        Switzerland\\
        email: schmitz@math.ethz.ch
}
\date{submitted 6 October 2005\\ accepted for publication in the 
``Electronic Journal of Probability'' \\10 July 2006}
\maketitle

\bibliographystyle{alpha}
\begin{abstract}
With the help of the methods developed in \cite{schmitz}, we highlight condition $(T)$ as a source of new examples of 'ballistic' diffusions in a random environment when $d \geq 2$ ('ballistic' means that a strong law of large numbers with non-vanishing limiting velocity holds). 
In particular we are able to treat the case of non-constant diffusion coefficients,
a feature that causes problems.
Further we
recover the ballistic character of two important classes of diffusions in a random environment 
by simply checking condition $(T)$. 
This not only points out to the broad range of examples where condition $(T)$ can 
be checked, but 
also fortifies our belief that
condition $(T)$ is a natural contender for the characterisation of ballistic diffusions in  a random environment when $d \geq 2$.
\end{abstract}
\pagenumbering{arabic}
\section{Introduction}

Diffusions in random environment in higher dimensions are still  poorly understood. An important tool in the study of diffusions in random environment has been the
``method of the environment viewed from the particle'', cf. \cite{olla94}, \cite{olla01} and the references therein, 
with the serious drawback that it only applies in very specific situations 
where one can find an (most often explicit) invariant measure for the process of 
the environment viewed from the particle. 
However new ideas recently emerged in the discrete framework of random walks in random environment 
that apply in the general setting, see  for instance \cite{szn04}, \cite{zeit} and the references therein. 
In particular the renewal-type arguments introduced in Sznitman-Zerner \cite{szn-zer}, 
and related developments centered around condition $(T)$ in Sznitman \cite{szn01}, \cite{szn02}, \cite{szn03}
 shed some light on walks with non-vanishing limiting velocity.
Let us as well point out the recent breakthrough in the diffusive regime, cf. Sznitman 
and Zeitouni \cite{szn-zeit}, where new methods falling outside the framework of 
these renewal techniques are developed. \\
In Shen \cite{shen} a renewal structure in the spirit of \cite{szn-zer}
was implemented in the continuous space-time setting. 
In our recent work \cite{schmitz}, we built up on these results, and 
showed that condition $(T)$ is also instrumental in the continuous setting:
 when $d \ge 2$, it implies
a strong law of large numbers with non-vanishing limiting velocity 
(which we refer to as {\it ballistic} behavior) and an invariance principle governing corrections to the law of large numbers.\\ 
This article mainly follows two objectives. 
We first highlight condition $(T)$ as a source of new examples of 
ballistic diffusions in random environment when $d \geq 2$.
Second, we
rederive the well-known ballistic character of a 
class of diffusions in random environment with a 
divergence-free drift field by simply checking condition $(T)$. 
We also check condition $(T)$ for an anisotropic gradient-type diffusion, the   ballistic nature 
of which was shown in Shen \cite{shen}. This not only points out to the broad range of examples 
where condition $(T)$ can be checked, but also fortifies our belief that condition $(T)$ is a natural contender for the characterisation of ballistic behavior.\\

Before describing our results in more details, let us recall the setting.\\
The {\it random environment} $\om$ is described by a probability space $(\Omega,\mathcal{A},\mathbb{P})$.  
We assume that there exists a group $\{t_x:x \in \R{d}\}$ of transformations on $\Omega$, jointly measurable in $x, \om$, which preserve the probability $\mathbb P$:
\begin{equation}
  \label{eq:stationarity}
  t_x \mathbb P =\mathbb P \,.
\end{equation}
On $(\Omega,\mathcal{A},\mathbb{P})$ we consider random variables
$a(\cdot)$ and $ b(\cdot)$ with respective values in the space of symmetric $d \times d$ matrices and $\R{d}$, and we define
\begin{equation} 
\label{eq:b-sigma}
a(x,\om) \df a(t_{x}\om), \quad b(x,\om) \df b(t_{x}\om).
\end{equation}
We shall assume that $a(\cdot,\om)$ is uniformly elliptic and bounded, and that $b(\cdot,\om)$ is bounded,
i.e. there are constants  $\nu \ge 1$, $\beta>0$, such that for all $\om \in \Omega$, all $x,y \in \R{d}$, 
\begin{equation}
\label{eq:a}
\tfrac{1}{\nu}|y|^{2}\leq y^t \cdot a(x,\om)\,y \leq \nu |y|^{2}, \qquad 
\left| b(x,\om) \right| \leq \beta,\,\\
\end{equation}
where $|\cdot|$ denotes the Euclidean norm for vectors resp. for square 
matrices, and $y^t$ stands for the transposed vector of $y$. We assume as well that
 $a(\cdot,\om)$, $b(\cdot,\om)$ are 
Lipschitz continuous, i.e. there is a constant $K>0$ such that for all $\om \in \Omega,~x,y  \in \R{d}$,
\begin{equation}
\label{eq:Lipschitz}
|a(x,\om)-a(y,\om)|+|b(x,\om)-b(y,\om)| \leq K|x-y|\,.
\end{equation} 
We will further assume \emph{finite range dependence}, i.e. there is an $R>0$
such that, when we define the $\sigma$-field $\mathcal{H}_{F} \df 
\sigma\{a(x,\cdot), b(x,\cdot):x \in F\}$, $F$ a Borel subset of $\R{d}$, then 
\begin{equation}
\label{eq:R-separation}
\mathcal{H}_{A} \text{ and }  \mathcal{H}_{B} \text{ are $\mathbb P$-independent when } d(A,B)>R \,,
\end{equation}
where $d(A,B)= \inf\{|x-y|: x \in A, y \in B \}$.
We define the differential operator attached to $a(\cdot, \om), b(\cdot, \om)$,
\begin{equation}
  \label{eq:diff-operator}
  \cL_\om =\tfrac  {1}{2}\sum_{i,j=1}^d a_{ij}(x,\om) \partial^2_{ij}+
          \sum_{i=1}^d b_i(x, \om)\partial_{i}\,,
\end{equation}
and, for $\all$, we denote with $P_{x,\om}$ the unique solution to the martingale problem for $\cL_\om$ started at $x \in \R{d}$, see Theorem 4.3 p.146 in \cite{bass2}.
The laws $ P_{x,\om}$  describe the
diffusion in the environment $\om$. We write $E_{x,\om}$ for the corresponding  expectation. We denote with $(X_t)_{t \geq 0}$ the canonical process on $C(\R{}_+, \R{d})$, 
with $(\mathcal F_t)_{t \geq 0}$ the canonical filtration,
and $\cF$ the Borel $\sigma$-field on $C(\R{}_+, \R{d})$.
The laws $ P_{x,\om}$ are usually called the {\it quenched} laws.
To restore some stationarity to the problem, it is convenient to introduce the \emph{annealed} laws $P_{x}$,  
which are defined as the semi-direct products: 
\begin{equation}
\label{eq:annealed}
P_{x} \df \mathbb{P} \times P_{x,\om},,~~ \mathrm{for}~ x \in \R{d}.
\end{equation}
The corresponding expectations are denoted with $E_{x}=\mathbb{E} \times E_{x,\om}$, 
where $\mE$ denotes the expectation with respect to the measure $\mP$. 
Observe that the Markov property is typically lost under the annealed laws. \\

Let us now explain in more details the purpose of this work. We first recall from \cite{schmitz} the definition of condition $(T)$ (where $T$ stands for transience).
Introduce, for $\ell$ a unit vector of $\R{d}$, $b, L>0$, the slabs 
$U_{\ell,b,L} \df \{x \in \R{d}:-bL < x \cdot \ell <L\}$. 
\begin{defn}
We say that \emph{condition $(T)$} holds relative to
$\ell \in S^{d-1}$, in shorthand notation $(T)|\ell$, if for all $\ell' \in S^{d-1}$ in a neighborhood of $\ell$, and for all $b>0$,
\begin{equation}
\label{eq:T}
\limsup_{L \to \infty} L^{-1} \log P_{0}[X_{T_{U_{\ell',b,L}}}\cdot \ell' < 0]<0\,,
\end{equation}
where $T_{U_{\ell',b,L}}$ denotes the exit time from the slab $U_{\ell',b,L}$.
\end{defn}
Let us recall that, when $d \geq 2$, and under our standing assumptions (\ref{eq:stationarity})-(\ref{eq:R-separation}), 
condition $(T)$ implies a ballistic strong law of large numbers and an annealed invariance principle governing corrections to the law of large numbers, see \cite{schmitz}:
\begin{equation}
\label{eq:lln}
P_{0}-\text{a.s.}, \quad \frac{X_{t}}{t} \rightarrow \tilde v, \quad \tilde v
\neq 0,  \text{ deterministic, with }  \tilde v \cdot \ell >0\,,  
\end{equation}
\begin{equation}
\label{eq:clt}
\begin{aligned}
&\text{and under $P_{0}$, $B_{\cdot}^{s}=\tfrac{X_{\cdot s }-\cdot s \tilde v}{\sqrt{s}}$
 converges in law on $C(\R{}_{+},\R{d})$, as $s \to \infty$, to a}\\
&\text{Brownian motion $B_{\cdot}$ with non-degenerate covariance matrix}.
\end{aligned}
\end{equation}

In Section 2, we provide a new class of examples where condition $(T)$ holds.
Namely, we prove that there is a constant $c_1(d,\nu,\beta,K,R)>1$ such that, for 
$\ell \in S^{d-1}$, the inequality
\begin{equation*}
\mE[(b(0,\omega) \cdot \ell)_{+}]>c_1~ \mE[(b(0,\omega) \cdot \ell)_{-}]
\end{equation*}
implies condition $(T)|\ell$, see Theorem \ref{thm:examples}. 
 This theorem extends the result of Theorem 5.2 in \cite{schmitz}
 in the sense that we drop the assumption  that the diffusion coefficient $a(\cdot,\om)$ 
equals the identity for all $\all$. This is more than a technical improvement, for 
this new set-up leads to severe difficulties.
Let us recall that in \cite{schmitz} we introduced for each bounded $C^\infty$-domain 
$U$ containing 0, an auxiliary diffusion with characteristics independent of the environment 
that exhibited the same exit distribution from $U$ as the annealed diffusion in the random 
environment
when starting at 0. With the objective of checking condition $(T)$, we approximated the slab 
$U_{\ell,b,L}$ by bounded $C^\infty$-domains, and using the auxiliary diffusion, 
we were able to restore some Markovian character to the question of controlling the
exit distribution from the slab $U_{\ell,b,L}$ under the annealed measure $P_0$. 
We keep this strategy in the present article, 
with the difference that now, 
when $d \ge 2$,
the construction of an auxiliary diffusion exhibiting the same properties as the 
auxiliary diffusion from \cite{schmitz} described 
above is 
more intricate.\\
Inspired by \cite{schmitz}, we 
are naturally led to choose for a bounded $C^\infty$-
domain $U$ containing 0
the following auxiliary diffusion and drift coefficients: 
\begin{equation}
\label{eq:a'}
\begin{aligned}
a'_U(x)&=\tfrac{\mE[g_{\om,U}(x)a(x,\om)]}{\mE[g_{\om,U}(x)]}\,,
  \text{if $x \in U \smallsetminus \{0\}$}, \quad a'_U(x)=Id, \text{ else},\\
b'_U(x)&=\tfrac{\mE[g_{\om,U}(x)b(x,\om)]}{\mE[g_{\om,U}(x)]}\,,
  \text{if $x \in U \smallsetminus \{0\}$}, \quad \,b'_U(x)=0, \,\,\,
\text{  else},
\end{aligned} 
\end{equation}
where $g_{\om,U}$ is the Green function corresponding to the quenched diffusion started 
at 0 and killed when exiting $U$. Notice that, when $a(\cdot,\om) \equiv Id$ for 
all $\all$ as in 
\cite{schmitz}, then $a'_U \equiv Id$, and
we recover the auxiliary diffusion from \cite{schmitz}, i.e. a Brownian motion perturbed
by a bounded measurable
drift $b'_U$. In the present setting, $a'_U$ is uniformly elliptic and bounded, but in general
 discontinuous at 0 and at the boundary $\partial U$. 
The massive discontinuity at the boundary may very well imply 
nonuniqueness for the martingale problem attached to $a'_U$ and $b'_U$.  \\ 
Thus, in the key Proposition \ref{prop:exit}, we prove the 
{\it existence} of a solution $P'_{0,U}$ to the martingale problem for $\cL'$ started 
at 0 such that
the respective exit distributions from $U$ under $P'_{0,U}$ and 
the annealed measure $P_0$ agree. Any such measure will serve as an auxiliary measure.\\
We cannot adapt the methods used in \cite{schmitz} to prove Proposition \ref{prop:exit} 
in the present setting. Due to the discontinuity of $a'_U$ at 0, a certain Dirichlet problem 
in the domain $U$, attached to the auxiliary coefficients $a'_U$, $b'_U$, may fail to 
have a strong solution in the Sobolev space $W^{2,p}(U)$, $p \ge 1$. This is in contrast to the situation in \cite{schmitz}, 
where $a'_U \equiv Id$ implied
the existence of a solution with good differentiability properties to the above Dirichlet problem.
As a result, when choosing arbitrary smooth and bounded approximations
$b_n$ of the auxiliary drift $b'_U$, we obtained 
convergence of the smooth solutions,
together with their derivatives, 
of the perturbed Dirichlet problems attached to $a_n \equiv Id$ and $b_n$.
This enabled us to conclude the proof 
of the result corresponding to Proposition \ref{prop:exit} in the setting 
of \cite{schmitz}, see also Remark \ref{rem:2.1}. \\
To avoid controls on derivatives in the present setting, we use results of Krylov \cite{krylov92}, \cite{krylov04}.
The specific choice of the auxiliary coefficients 
$a'_U,\,b'_U$ enables us to construct a diffusion started at 0 and attached to $a'_U,\,b'_U$,
that, 
killed when exiting $U$, admits 
$\mE[g_{\om,U}]$
as occupation time density in $U$. This implies that the above diffusion exhibits the 
same exit distribution from $U$ as the annealed diffusion in random environment.\\


In Section 3, we study a Brownian motion perturbed by a random divergence-free drift field
the $\mathbb P$-expectation of which does not vanish.
The ballistic character of this class of diffusions (including time-dependent 
drift fields) follows for instance from the results in Landim, Olla, Yau \cite{landim-olla-yau} 
via the method of ``the environment viewed from the particle''.
We use the pointwise Gaussian controls on the quenched semigroup from Norris 
\cite{norris} (see also Osada \cite{osada1} for estimates that only accommodate a drift 
term with vanishing $\mathbb P$-expectation), that enable us to 
 check condition $(T)$ by a standard Markov-process-type computation, see Theorem \ref{thm:div}. 
With the help of our results in \cite{schmitz}, we thus rederive
the 
 ballistic behavior (in the time-independent setting) of this class of diffusions 
in random environment. Let us also mention that, in the case where the $\mathbb P$-expectation
of the drift is zero, diffusive behavior has been shown, see for instance \cite{fan-papa},\cite{oel}, \cite{osada}
for the time-independent setting, and \cite{fan-kom}, \cite{kom-olla-01}, \cite{landim-olla-yau} for the time-dependent setting.

In Section 4, we consider Brownian motion perturbed by a random drift that can be written as the gradient of a {\it non-stationary} scalar potential. 
The ballistic behavior of this class of diffusions was only recently shown in Shen \cite{shen}. Whereas 
the method of ``the environment viewed from the particle'' applies successfully when the potential is {\it stationary}, see \cite{olla94} 
or section 3.7 in \cite{olla01}, it does not apply in the above setting.
 Shen shows the existence of the first two moments of a certain regeneration time, see Theorem 4.11 in \cite{shen}, and obtains the ballistic character as a direct consequence of his previous results, see Theorems 3.1 and 3.2. 
of \cite{shen}. 
With closely related techniques, we 
are able to prove in Theorem \ref{thm:grad} that condition $(T)$ holds in this setting.\\

Let us also mention that in the context of the examples handled in Sections 3 and 4, 
the analogue of (\ref{eq:T}),
where the annealed measure $P_0$ is replaced by the quenched measure $P_{0,\om}$, holds
uniformly in $\om$, cf. (\ref{eq:exit-0}) and (\ref{eq:exit}).\\

{\bf Convention on constants:}
Unless otherwise stated, constants only depend on the quantities $\beta, \nu ,K, R, d $. 
Dependence on additional parameters appears in the notation. 
Generic positive constants are denoted by $c$.
When constants are not numerated, their value may change from line to line.\\

{\bf Acknowledgement:}
I am sincerely indebted to my advisor Prof. A.-S. Sznitman for his constant advice during the completion of this work. 
I would like to thank as well Laurent Goergen for fruitful discussions, and the referee for helpful comments.

\section{A new class of examples where Condition $(T)$ holds}
\label{sec:2}
We first introduce some additional notation. For $z \in \R{}$, we denote with $\lfloor z \rfloor \df \sup\{k \in \mathbb Z:k \le z\}$ the rounded value of $z$. 
We denote with  $|\cdot|$ the Euclidean norm, and with $|\cdot|_\infty$ the supremum norm.
For $x \in \R{d}$, $r>0$, we write 
$B_{r}(x)$ for the open Euclidean ball with radius $r$ and center $x$. 
We call an open connected subset of $\R{d}$ a domain, and for a subset $A$  
of $\R{d}$, we write
 $A^c$ for its complement, $\bar A$ for its closure, and $\partial A$ for its boundary. The 
diameter of $A$ is defined as $diam\,A=\sup\{d(x,y):x,y \in A\}$.
For two sets $A,B$ in $\R{d}$, $d(A,B)= \inf\{|x-y|: x \in A, y \in B \}$ stands 
for the distance between $A$ and $B$.
For an open subset $U$ of $\R{d}$, the $(\cF_t)$-stopping time $T_U= \inf{\{t \geq 0:X_t \notin U\}}$ denotes the exit time from $U$.  \\
Given measurable functions $a,\,b$ on $\R{d}$ with values in the space of 
symmetric $d \times d$ matrices, resp. in $\R{d}$, 
that satisfy the conditions in (\ref{eq:a}) with the respective constants $\nu, 
\beta$ (with $a(\cdot),\,b(\cdot)$ in place of $a(\cdot,\om),\,b(\cdot,\om)$), we say that the attached differential operator
$\cL$ is of class $\cN(\nu,\beta)$.\\
Recall the convention on constants at the end of the introduction. The main result of this section is
\begin{thm}
\label{thm:examples}
($d \geq 1$) \\
Assume (\ref{eq:stationarity})-(\ref{eq:R-separation}). There is a constant $c_1>1$, such that for $\ell \in S^{d-1}$, the 
inequality
\begin{equation}
\label{eq:criterion}
\mE[(b(0,\omega) \cdot \ell)_{+}]>c_1~ \mE[(b(0,\omega) \cdot \ell)_{-}]
\end{equation}
implies $(T)|\ell$.
\end{thm}
This is an extension of Theorem 5.2 in our previous work \cite{schmitz}, 
where we proved the same statement under the additional assumption that 
$a \equiv  Id$.
\begin{rem} \rm
One can easily argue that $c_1 \ge 1$, and the proof of Theorem  \ref{thm:examples} shows in fact 
that $c_1>1$. The true value of the constant is unknown. 
However, when $\mathbb P$-a.s., $(b(0,\om)\cdot \ell)_-= 0$, then Theorem \ref{thm:examples} immediately implies,
regardless of the value of $c_1$:
\begin{equation}
\label{eq:non-nestling}
\begin{aligned}
&\text{Condition $(T)$ holds when $d \geq 1$ and there is $\ell \in S^{d-1}$ and $\delta >0$,} \\
&\text{such that $b(0,\om)\cdot \ell \geq 0$ for all $\om \in \Omega$, and  
 $p_\delta = \mathbb P[b(0,\om)\cdot \ell \geq \delta]>0$}\,.
\end{aligned}
\end{equation}
If there is $\delta >0$ such that $p_\delta=1$, this is in the spirit of the 
{\it non-nestling} case, which is in fact already covered by Proposition
5.1 in \cite{schmitz}, and else, of the {\it marginal nestling} case
in the discrete setting,  see Sznitman \cite{szn00}.
Of course, Theorem \ref{thm:examples} also comprises more involved examples of condition $(T)$ where 
$b(0,\om)\cdot \ell$ takes both positive and negative values for every $\ell \in S^{d-1}$. They  correspond to the {\it plain-nestling} case in \cite{szn00}.
When $d \geq 2$, and when the covariance matrix $a$ is not the identity, the examples of condition $(T)$ corresponding
to the marginal- resp.~plain-nestling case provide new examples of ballistic diffusions in random environment.
$\hfill \Box$
\end{rem}
The main step in the proof of Theorem \ref{thm:examples} is to construct for each 
bounded $C^\infty$-domain $U$ containing 0 an 
auxiliary diffusion with characteristics independent of the environment, that, when 
started at 0, exhibits the same exit distribution from $U$ as the annealed diffusion 
in random environment started at 0. We provide this key result in subsection 2.1, 
cf. Proposition \ref{prop:exit}. \\
In subsection \ref{subsec:T}, we prove Theorem \ref{thm:examples} using slight variations
on the methods from \cite{schmitz}.
The auxiliary diffusion enables us to restore some Markovian character to the task
of checking condition $(T)$.
First, we show that condition $(T)$ is implied by a certain condition $(K)$, 
which is similar to Kalikow's condition in the discrete set-up.
We then check condition $(K)$, which finishes the proof of Theorem \ref{thm:examples}.

\subsection{An auxiliary diffusion and its exit distribution}
Unless otherwise stated, $U$ denotes from now on a bounded $C^\infty$-domain containing 
0.
Recall the subtransition density  $p_{\om,U}(s,x,\cdot)$ at 
time $s > 0$ for the quenched diffusion started at $x \in \R{d}$, that fulfills for each open 
set $G \subset U$, 
$P_{x,\om}[X_s \in G, T_U >s]=\int_G p_{\om,U}(s,x,y)\,dy$.
It can for instance be defined via Duhamel's formula, see \cite{str} p.331:
\begin{equation}
\label{eq:duhamel}
p_{\om,U}(s,x,y)=p_\om(s,x,y)-E_{x,\om}[p_\om(s-T_U,X_{T_U},y), T_U < s],
\quad s>0,\,\,x,y \in \R{d},
\end{equation}
where $p_\om(s,x,\cdot)$ denotes the transition density corresponding
to the quenched diffusion started at $x \in \R{d}$. 
One then obtains the Green function for the quenched diffusion started at $0$ and killed when exiting $U$ via
\begin{equation}
  \label{eq:green-function}
  g_{\om,U}(x) = \int_{0}^{\infty} p_{\omega, U}(s,0,x)\,ds, \quad x \in \R{d}.
\end{equation}
We define the auxiliary diffusion and drift coefficients, see also (\ref{eq:a'}), through
\begin{equation}
\label{eq:b'}
\begin{aligned}
a'_U(x)&=\tfrac{\mE[g_{\om,U}(x)a(x,\om)]}{\mE[g_{\om,U}(x)]}\,,
  \text{if $x \in U \smallsetminus \{0\}$}, \quad a'_U(x)=Id, \text{ else},\\
b'_U(x)&=\tfrac{\mE[g_{\om,U}(x)b(x,\om)]}{\mE[g_{\om,U}(x)]}\,,
  \text{if $x \in U \smallsetminus \{0\}$}, \quad \,b'_U(x)=0, \,\,\,
\text{  else},
\end{aligned}
\end{equation} 
and we denote with $\cL'$ the attached differential operator.\\
Notice that the definition of $a'_U,\,b'_U$ in 0 and outside $U$ is quite arbitrary.
Since we are interested in the exit distribution from $U$ of a suitable diffusion started 
at 0 and
attached to $a'_U,b'_U$, cf. Proposition \ref{prop:exit} below, we only 
require that $a'_U,\,b'_U$ are uniformly elliptic and bounded outside $U$.
Let us further mention that changing the value of $a'_U(0)$ within the class of symmetric and elliptic matrices does not affect issues such 
as existence and uniqueness of the solution to the martingale problem for $\mathcal{L}'$, see for instance \cite{bass2} p. 149 and 
Theorem 1.2 p.132 therein.\\
The next lemma will be useful in the remainder of this section.
\begin{lemma}
\label{lemma:L'}
$\cL' \in \mathcal N(\nu,\beta)$. 
Moreover,  $a'_U$, $b'_U$, and $g_{\om,U}$ are continuous at  $U \smallsetminus\{0\}$.
\end{lemma} 
\begin{proof}
It follows from (\ref{eq:a}) and from the definition of $a'_U, b'_U$ 
that $\cL' \in \mathcal N(\nu, \beta)$. 
Theorem 4.5 p.141 in \cite{fried} states that 
the transition density $p_\om(s,x,y)$ is jointly continuous in $s>0,\,x,y \in \R{d}$. 
By Duhamel's formula, see (\ref{eq:duhamel}), and dominated convergence, 
it follows that the subtransition $p_{\om,U}(s,x,y)$ is continuous in $y
\in U$.
From Lemma \ref{lemma:PDE} in the Appendix, and from similar computations as carried out between (\ref{eq:sup-green}) and (\ref{eq:chap-kolm}) in the 
Appendix, and applying dominated convergence, we obtain the continuity of $g_{\om,U}$ in $U \smallsetminus \{0\}$. Further, 
by the upper bound in Corollary \ref{cor:green}  
and dominated 
convergence, and
 by continuity of $a(\cdot,\om)$ resp. $b(\cdot,\om)$, see (\ref{eq:Lipschitz}), 
 we see that
$a'_U$ resp. $b'_U$ are continuous at  $U \smallsetminus \{0\}$.
\end{proof}
The following proposition restores some Markovian character to the exit problem of
the diffusion $X_\cdot$ under the annealed measure $P_0$ (recall $P_0$ in 
(\ref{eq:annealed})).
\begin{prop}
\label{prop:exit}($d \ge 1$)
Let $U$ be a bounded $C^{\infty}$-domain containing 0. Then there is a solution $P'_{0,U}$
to
 the martingale problem for $\cL'$ started at 0, 
such that $X_{T_U}$ has same law under $P'_{0,U}$ 
and the annealed measure $P_0$. 
\end{prop}
In the sequel we will denote with $P'_{0,U}$ any such measure and refer to it as 
an auxiliary measure, and we will 
denote with $E'_{0,U}$ the corresponding expectation.
\begin{rem} \rm
\label{rem:1}
When $d \le 2$, Exercises 7.3.3 and 7.3.4 p.192 in \cite{str-var}, and the fact that
$\cL' \in \cN(\nu, \beta)$,  
show that the martingale problem for $\cL'$
is well-posed, so that the auxiliary 
measure $P'_{0,U}$ corresponds to the unique solution to the martingale problem. 
In higher dimensions, uniqueness to the martingale problem may fail, see Nadirashvili \cite{nad96}.
Uniqueness holds for instance when the 
set of discontinuities of the diffusion coefficient
contains only a finite number of cluster points, see Krylov \cite{krylov04}.
In our case, where the diffusion coefficient $a'_U$ is in general discontinuous
at the boundary $\partial U$, the question of uniqueness is open. 
However, 
since $a'_U$ 
is continuous at $U \smallsetminus \{0\}$, 
one can in fact show with the above result from Krylov \cite{krylov04} that 
 two arbitrary solutions to the martingale problem for $\cL'$ 
started at 0 agree on $\cF_{T_U}$, although we will not prove it here. 
As a consequence, every solution to the martingale problem for $\cL'$ started at 
0 can be used as auxiliary measure, but we will not use this here.
$\hfill{}\Box$
\end{rem}

\subsubsection{The case $d=1$}
\begin{proof}[{\bf Proof of Proposition \ref{prop:exit} when $d=1$}]
According to Remark \ref{rem:1}, we denote with $P'_{0,U}$ the unique solution to
the martingale problem for $\cL'$ started at $0$. 
Corollary 4.8 p.317 in \cite{kar-shr}, combined with Remark 4.3 p.173
therein, shows the existence of a Brownian motion $W$ defined on
$(C(\R{}_+,\R{}), \cF, P'_{0,U})$ such that
\begin{equation}
\label{eq:35}
P'_{0,U}-\text{a.s.}, \quad X_t=\int_0^t \sqrt{a'_U(X_s)}\,dW_s
+\int_0^t b'_U(X_s)\,ds.
\end{equation}
Say that $U=(\alpha_1,\alpha_2)$, $\alpha_1<0<\alpha_2$, and define a scale function $s$, 
\begin{equation}
\label{eq:31}
s(x)=\int_{\alpha_1}^x \exp\{-2\int_{\alpha_1}^y \frac{b'_U(z)}{a'_U(z)}\,dz\}\,dy,\quad 
x \in \R{}.
\end{equation}
Notice that $s \in C^1(\R{})$, and that, due to Lemma \ref{lemma:L'}, the second derivative
of $s$ exists outside $A=\{\alpha_1,0,\alpha_2\}$, and satisfies 
\begin{equation}
\label{eq:mart0}
\partial^2_x s(x)=-\frac{2\,b'_U(x)}{a'_U(x)}\,\partial_x s(x),\quad x \in A^c,\quad 
\text{and hence,}\quad\cL's=0 \text{ on }A^c.
\end{equation}
On $A$, we define $\partial^2_x s$ through the equality in (\ref{eq:mart0}),
 so that $\cL' s=0$ on $\R{}$.
Applying the generalised Ito rule, cf. problem 7.3 p.219 in \cite{kar-shr},
 and taking expectations, we find that, 
\begin{equation}
\label{eq:mart1}
E'_{0,U}[s(X_{t \wedge T_U})]=s(0)+E'_{0,U}
[\int_0^{t \wedge T_U}\cL's(X_u)\,du]=s(0),\quad t \ge 0.
\end{equation}
Using $E'_{0,U}[T_U]<\infty$, together 
with dominated convergence, we obtain
\begin{equation}
\label{eq:mart3}
E'_{0,U}[s(X_{T_U})]=s(0).
\end{equation}
A similar equation as (\ref{eq:35}) holds with $P_{0,\om}$ in place of
$P'_{0,U}$, and $a(\cdot,\om),\,b(\cdot,\om)$ in place of $a'_U(\cdot),\,b'_U(\cdot)$ respectively.
We proceed similarly as in (\ref{eq:mart1}) and below, and obtain
from the definition of the Green function $g_{\om,U}$, cf. (\ref{eq:green-function}),
\begin{equation}
\label{eq:mart4}
E_{0,\om}[s(X_{T_U})]=s(0)+E_{0,\om}
[\int_0^{T_U}\cL_\om s(X_u)\,du]\stackrel{(\ref{eq:green-function})}{=}
\int_{\alpha_1}^{\alpha_2} g_{\om,U}(x)\,\cL_\om s(x)\,dx.
\end{equation}
We integrate out (\ref{eq:mart4}) with respect to $\mathbb P$, and find from the definition 
of $a'_U,\,b'_U$, see (\ref{eq:b'}),
\begin{equation}
\label{eq:mart5}
E_0[s(X_{T_U})]=s(0)+\int_{\alpha_1}^{\alpha_2} \mE[g_{\om,U}(x)]\,\cL' s(x)\,dx=s(0).
\end{equation}
Since $s(\alpha_1)=0$, (\ref{eq:mart3}) and (\ref{eq:mart5}) imply that
$P'_{0,U}[X_{T_U}=\alpha_2]=P_0[X_{T_U}=\alpha_2]$, which is our claim.
\end{proof}

\subsubsection{The case $d \ge 2$}
The higher-dimensional case is more intricate. 
In Remark \ref{rem:discontinuous}, we show that when $d \ge 3$, then $a'_U$ is in general discontinuous at the origin.
In Remark \ref{rem:2.1}, we first recall from \cite{schmitz} the main points of the proof of 
Proposition \ref{prop:exit} when in addition $a= Id$ holds.
Then we explain that these methods fail when $a'_U$ is discontinuous. Hence, 
in the setting of a general covariance matrix $a$, where $a'_U$ is in general discontinuous,
the arguments from \cite{schmitz} break down. 
\begin{rem} \rm
\label{rem:discontinuous}
Theorem 4.1 p.80 in \cite{boboc-mustata} 
   states that there is a constant $c >0$, such that for all $\all$, $\varepsilon >0$, there is $\eta >0$ such that if $x \in U$, $|x| \le \eta$, then
   \begin{equation}
     \label{eq:boboc}
     1-\varepsilon \le \frac{g_{\om,U}(x)}{c\, h_\om(x)} \le 1+\varepsilon\,,
   \end{equation}
   where $ h_\om(x)=(x^t \cdot a(0,\om)\,x)^{\frac{2-d}{2}}$.
  Choose a real-valued sequence $\{r_n\}_n$ converging to 0,
   and define sequences $x_n^{(i)}=r_n e_i$, where $e_i$ is the $i$-th unit vector of 
   the canonical 
   basis. We find with the help of (\ref{eq:boboc}) and dominated convergence that
   \begin{equation*}
   \lim_n a'_U(x_n^{(i)})= \mE[a_{ii}(0,\om)^{\frac{2-d}{2}}a(0,\om)]\Big /
   \mE[a_{ii}(0,\om)^{\frac{2-d}{2}}]\,, 
   \end{equation*}
  so that the limit depends on $i \in \{1,\ldots,d\}$. Hence $a'_U$ is in general discontinuous at  0. \\
Let us point out that there are examples of diffusions, attached to uniformly elliptic, 
but discontinuous diffusion coefficients, that return to their starting point
with probability one, cf. Proposition 3.1 p.104 in \cite{bass2}.
$\hfill{}\Box$
\end{rem}

\begin{rem}\rm
\label{rem:2.1}
1. Let us recall that in \cite{schmitz}, the additional 
assumption that $a \equiv Id$ in Theorem 5.2 allowed us 
to choose for each bounded $C^\infty$-domain $U$ containing 0 a Brownian motion perturbed 
by a bounded measurable drift $b'_U$ (defined as in (\ref{eq:b'})) as the (uniquely 
defined) auxiliary diffusion. We emphasize that in this setting $a'_U=Id$.\\ 
We recall the main points in the proof of Proposition 5.4 in \cite{schmitz} (which corresponds to Proposition \ref{prop:exit} when $a = Id$): 
from the martingale problem, and by definition of $b'_U$, we derive a sort 
of ``forward'' equation, namely for all $f \in C^2(\bar U)$,
\begin{equation}
\label{eq:forwardeq}
E_0[f(X_{T_U})]-E'_{0,U}[f(X_{T_U})]= 
\int_U (\mE[g_{\om,U}(x)]-g'_U(x))(\tfrac{1}{2} \,\Delta+ b'_{U}(x)\,
\nabla )f(x)\,dx\,,
\end{equation}
where $g_{\om,U}$ and $g'_U$ are the respective Green functions corresponding to 
the quenched resp. auxiliary diffusion started at 0 and killed when exiting $U$.
We choose an arbitrary smooth function $\phi$, 
and we denote with $u$ the unique strong solution in the Sobolev space 
$W^{2,p}(U)$, $p>d$, to 
the problem 
\begin{equation}
\label{eq:dir1}
1/2 \,\Delta\, u + b'_{U}\,\nabla u=0 \text{ in }U, \quad u=\phi \text{ on }\partial U.
\end{equation} 
By considering the perturbed Dirichlet problems 
$1/2 \,\Delta\, u_n + b'_{U,n}\,
\nabla u_n=0$ in $U$, $u_n=\phi$ on $\partial U$,  where the  $b'_{U,n}$ 
are smooth and converge a.e. in $U$ to 
the bounded measurable drift $b'_U$, we obtain smooth 
solutions $u_n$ that converge in $W^{2,p}(U)$ to $u$.
Sobolev's inequality yields an upper bound, uniform in $n \ge 1$,
of the supremum norm on $U$ of the gradients of the $u_n$.
Since on $U$,  $1/2 \,\Delta\, u_n + b'_U\,\nabla u_n=
(b'_U\,-b'_{U,n}\,)\nabla u_n$ holds, we then conclude the proof by
inserting the smooth functions $u_n$ in (\ref{eq:forwardeq}) and  using dominated convergence.\\

2. Unlike the above setting (where $a'_U=Id$), $a'_U$ is now in general discontinuous at the origin when $d \ge 3$, see Remark \ref{rem:discontinuous}.
However,  the existence of a unique
strong solution in the Sobolev space $W^{2,p}(U)$, $p \ge 1$, 
for the Dirichlet problem corresponding to (\ref{eq:dir1}),
\begin{equation}
\label{eq:dir2}
\cL' u=0 \text{ in }U, \quad u=\phi \text{ on }\partial U,
\end{equation} 
is only guaranteed when  $a'_U$ is continuous on $\bar U$, see \cite{gil-tru} p.241.
Hence, we have no controls on derivatives as above, and the arguments from \cite{schmitz} break down.
$\hfill{}\Box$ 
\end{rem}
To avoid controls on derivatives, we use a result of Krylov \cite{krylov04}.
The specific choice of the auxiliary coefficients 
$a'_U,\,b'_U$ enables us to 
find a diffusion started at 0 and attached to $a'_U,\,b'_U$ that, killed when exiting $U$, admits 
$\mE[g_{\om,U}]$
as occupation time density in $U$, and we show that this diffusion exhibits the 
same exit distribution from $U$ as the annealed diffusion in random environment.
\begin{proof}[{\bf Proof of Proposition \ref{prop:exit} when $d \ge 2$}]
From the martingale problem, we find for $f \in C^2(\bar U)$:
\begin{equation}
\label{eq:2.7}
E_{0,\om}[f(X_{t \wedge T_U})]=f(0)+E_{0,\om}[\int_0^{t \wedge T_U}\cL_\om f(X_s)\,ds].
\end{equation}
Using ellipticity and standard martingale controls, it is immediate that $\sup_\om E_{0,\om}[T_U]<\infty$. By dominated convergence, and 
by the definition 
of the quenched Green function $g_{\om,U}$, cf. (\ref{eq:green-function}), 
we find that
\begin{equation}
\label{eq:2.8}
E_{0,\om}[f(X_{T_U})]=f(0)+\int_U g_{\om,U}(y)\cL_\om f(y)\,dy.
\end{equation}
After $\mP$-integration, and from the definition of $\cL'$, cf. (\ref{eq:b'}), we 
obtain that
\begin{equation}
\label{eq:2.9}
E_0[f(X_{T_U})]=f(0)+\int_U \mE[g_{\om,U}(y)]\cL'f(y)\,dy.
\end{equation}
In particular the assumption of Theorem 2.14 in \cite{krylov04} is satisfied (with 
the respective choice of the domains $Q=\R{d}$ and $D=U$ in the notations of \cite{krylov04}),
 and we infer the existence of a process 
$X$ and a Brownian motion $W$, defined on some probability space
$(C,\mathcal C,Q_0)$,  such that 
\begin{equation}
\label{eq:2.11}
Q_0-\text{a.s.}, \quad  X_t =\int_0^t \sigma'_U (X_s)\,dW_s+
 \int_0^t b'_U(X_s)\,ds\,,
\end{equation}
where $\sigma'_U$ is a measurable square root of $a'_U$, i.e. 
$a'_U=\sigma'_U \sigma'^t_U$, and such that for every bounded measurable
function $\psi$,
\begin{equation}
\label{eq:repres0}
\int_U \mE[g_{\om,U}(y)]\psi(y)\,dy = E'_{0,U}[\int_0^{T_U}\psi(X_s)\,ds],
\end{equation}  
where $E'_{0,U}$ denotes the expectation with respect to the measure
$P'_{0,U}$ induced by $X$ on $(C(\R{}_+,\R{d}),\cF)$.
Corollary 4.8 p.317 in \cite{kar-shr} shows that
$P'_{0,U}$  solves the martingale problem for $\cL'$ started at $0$. 
Insert $\psi =1$ in (\ref{eq:repres0}) to find $E'_{0,U}[T_U]=E_0[T_U]<\infty$, see 
below (\ref{eq:2.7}). 
For all $f \in C^2(\bar U)$, we derive a similar equation as (\ref{eq:2.7}) under the 
measure $P'_{0,U}$, 
so that, when letting $t \to \infty$, and using dominated convergence, 
we find that
\begin{multline*}
E'_{0,U}[f(X_{T_U})]=f(0)+E'_{0,U}[\int_0^{T_U}\cL'f(X_s)\,ds]\\
\stackrel{(\ref{eq:repres0})}{=}f(0)+\int_U \mE[g_{\om,U}(y)]\cL'f(y)\,dy 
\stackrel{(\ref{eq:2.9})}{=}E_0[f(X_{T_U})].
\end{multline*}
Since $f \in C^2(\bar U)$ is arbitrary, the claim of the proposition follows.
\end{proof}
\begin{rem} \rm
Although we will not use it here, one can show that in the terminology of 
Krylov \cite{krylov92}, \cite{krylov04}, $\mE[g_{\om,U}]$ is the unique
(up to Lebesgue equivalence) Green function of $\cL'$ in $U$ with pole at 
0. (This is a rather straightforward consequence of (\ref{eq:repres0}) and the fact 
that any two solutions to the martingale problem for $\cL'$ started at 0 agree on 
$\cF_{T_U}$, cf. Remark \ref{rem:1}).
$\hfill{}\Box$
\end{rem}

\subsection{The proof of Theorem \ref{thm:examples}}
\label{subsec:T}
We now recall from \cite{schmitz} the definition of condition $(K)$, which is similar to Kalikow's condition in the discrete set-up. 
With the help of Proposition {\ref{prop:exit}}, we show that it implies condition $(T)$. The proof 
of Theorem \ref{thm:examples} is finally concluded by checking condition $(K)$. 
The proofs in this subsection are adaptations from the corresponding proofs
in section 5 in \cite{schmitz}.
\begin{defn}
\label{def:K}
($d \ge 1$) Let $\ell \in S^{d-1}$. We say that condition ($K)|\ell$ holds, if there is an
$\epsilon >0$, such that for all bounded domains $U$ containing 0 
\begin{equation}
\label{eq:K}
 \inf_{x \in U \smallsetminus \{0\}, \text{ d} (x,\partial U)>5R} b'_U (x)\cdot \ell > \epsilon \,,
\end{equation}
with the convention $\inf \varnothing = +\infty$.
\end{defn}
As alluded to above, the next step is
\begin{prop}
\label{prop:K}
$(K)|\ell \Rightarrow (T)|\ell$\,. 
\end{prop}

\begin{proof}
The set of $\ell \in S^{d-1}$ for which (\ref{eq:K}) holds is open and  hence our claim will follow if for such an $\ell$ we show that
\begin{equation}
  \label{eq:T|l}
  \limsup_{L \to \infty}L^{-1}\log P_0[X_{T_{U_{\ell,b,L}}}\cdot \ell <0]<0\,.
\end{equation}
Denote with $\Pi_\ell(w) \df w-(w\cdot \ell) \ell$, $w \in \R{d}$, the projection on the orthogonal complement of $\ell$, and define 
\begin{equation}
  \label{eq:bounded-set}
  V_{\ell,b,L} \df \left\{x \in \R{d}: -bL < x \cdot \ell < L, |\Pi_\ell(x)| < L^{2} \right\}\,.
\end{equation}
In view of Proposition \ref{prop:exit}, we choose bounded 
$C^\infty$-domains $\tilde  V_{\ell,b,L}$ such that 
\begin{equation}
  \label{eq:set-inclusion}
   V_{\ell,b,L} \subset \left\{x \in \R{d}: -bL < x \cdot \ell < L, |\Pi_\ell(x)| < L^{2} +5R \right\} \subset \tilde  V_{\ell,b,L} \subset U_{\ell,b,L}\,.
\end{equation}
(When $d=1$, $\Pi_\ell(w) \equiv 0$, and we simply have that $U_{\ell,b,L}=V_{\ell,b,L}=\tilde V_{\ell,b,L}$.)
We denote with $P'_{0,\tilde V_{\ell,b,L}}$ an auxiliary measure, i.e. $P'_{0,\tilde V_{\ell,b,L}}$ 
solves the martingale problem for $\cL'$ started at 0, and $X_{T_{\tilde V_{\ell,b,L}}}$ 
has same law under $P'_{0,\tilde V_{\ell,b,L}}$ and $P_0$, cf. Proposition \ref{prop:exit}.
To prove (\ref{eq:T|l}), it will suffice to prove that
\begin{equation}
  \label{eq:exit-bounded-set}
  \limsup_{L \to \infty}L^{-1}\log  P'_{0,\tilde V_{\ell,b,L}}[X_{T_{V_{\ell,b,L}}}
  \cdot \ell <L]<0\,.
\end{equation}
Indeed, once this is proved, it follows from (\ref{eq:set-inclusion}) that
\begin{equation}
  \label{eq:exit-bounded-set2}
  \limsup_{L \to \infty}L^{-1}\log  P'_{0,\tilde V_{\ell,b,L}}[X_{T_{\tilde V_{\ell,b,L}}}\cdot \ell <L]<0\,.
\end{equation}
Then, by construction of $ P'_{0,\tilde V_{\ell,b,L}}$, (\ref{eq:exit-bounded-set2})
holds with $ P'_{0,\tilde V_{\ell,b,L}}$ replaced by $P_0$, and, using (\ref{eq:set-inclusion}) once
more, (\ref{eq:T|l}) follows.\\
We now prove (\ref{eq:exit-bounded-set}). 
By (\ref{eq:set-inclusion}) and (\ref{eq:K}), we see that for $x \in V_{\ell,b,L}$,
\begin{equation}
  \label{eq:b'-bounds}
  b'_{\tilde V_{\ell,b,L}}(x)\cdot \ell \geq 
\begin{cases}
 \epsilon,  &\text{ if } -bL+5R<x \cdot \ell <L-5R \text{ and } x \neq 0,\\ 
 -\beta,   &\text{ else }.
\end{cases}
\end{equation}
We thus consider the  process $X_t \cdot \ell$. 
We introduce the function $u(\cdot)$ on $\R{}$, which is
defined on $[-bL,L]$ through
\begin{equation}
u(r) \df
\begin{cases}
  \alpha_1 e^{\alpha_2 \tfrac{\epsilon}{\nu} (bL-5R)}(\alpha_3-e^{4\nu \beta (r-(-bL+5R))}),
  &\text{if } r \in [-bL,-bL+5R]\,,\\
  e^{-\alpha_2 \tfrac{\epsilon}{\nu} r}, &\text{if } r \in (-bL+5R,L-5R)\,,\\
  \alpha_4 e^{-\alpha_2 \tfrac{\epsilon}{\nu} (L-5R)}(\alpha_5-e^{4\nu \beta (r-(L-5R))}),
  &\text{if } r \in [L-5R,L]\,,
\end{cases}  
\end{equation}
and which is extended boundedly and in a $C^2$ fashion outside $[-bL,L]$, and such that $u$ is twice differentiable in the points $-bL$ and $L$.
The numbers $\alpha_i$, $1\leq i \leq 5$, are chosen positive and 
independent of $L$, via 
\begin{equation}
  \label{eq:coefficients}
  \alpha_5=1+e^{20\nu \beta R},\, \alpha_4=e^{-20\nu \beta R},\, \alpha_2 = \min (1, \frac{4 \nu^2 \beta}{\epsilon}e^{-20\nu \beta R}),\,
 \alpha_1=\frac{\epsilon \alpha_2}{4 \nu^2 \beta},\, \alpha_3=1+\frac{4 \nu^2\beta}{\epsilon \alpha_2}\,.
\end{equation}
Then, on $[-bL,L]$, $u$ is positive, continuous and decreasing.
In addition, one has with the definition $j(r)=u'(r_+)-u'(r_-)$, 
\begin{equation}
\label{eq:derivative}
j(-bL+5R)=0, \text{ and } j(L-5R) \le 0\,.
\end{equation}
On $\R{d}$ we define the function $\tilde u(x)=u(x \cdot \ell)$, and
for $\lambda$ real, we define on $\R{}_+ \times \R{}$ the function
$v_\lambda(t,r) \df e^{\lambda t}u(r)$, and on $\R{}_+ \times \R{d}$ the function
$\tilde v_\lambda(t,x) \df v_\lambda(t,x \cdot \ell)= e^{\lambda t}\tilde u(x)$. 
We will now find $\lambda_0$ positive such that
\begin{equation}
  \label{eq:supermartingale}
  v_{\lambda_0}(t \wedge T_{V_{\ell,b,L}}, X_{t \wedge T_{V_{\ell,b,L}}} \cdot \ell) \text{ is a positive supermartingale under }P'_{0,\tilde V_{\ell,b,L}}.
\end{equation}
Corollary 4.8 p.317 in \cite{kar-shr}, combined with Remark 4.3 p.173 therein, shows the existence of a $d$-dimensional Brownian motion $W$ defined on 
$(C(\R{}_+,\R{d}), \mathcal F,P'_{0,\tilde V_{\ell,b,L}})$, such that
\begin{equation*}
 P'_{0,\tilde V_{\ell,b,L}}-\text{a.s.}, \quad Y_t \df X_t \cdot \ell=
 \int_0^t \sigma'_{\tilde V_{\ell,b,L}}(X_s)\cdot \ell\,dW_s+
 \int_0^t b'_{\tilde V_{\ell,b,L}}(X_s)\cdot \ell\,ds\,,
\end{equation*}
where $\sigma'_{\tilde V_{\ell,b,L}}$ is a measurable square root of $a'_{\tilde V_{\ell,b,L}}$, i.e. 
$a'_{\tilde V_{\ell,b,L}}=\sigma'_{\tilde V_{\ell,b,L}}\sigma'^t_{\tilde V_{\ell,b,L}}$.
Writing $u$ as a linear combination of convex functions, we find from the 
generalised It\^o rule, see \cite{kar-shr} p.218, that
 \begin{align}
  \label{eq:ito}
P'_{0,\tilde V_{\ell,b,L}}-\text{a.s.,}\quad u(Y_{t})\,=\,1+\int_0^{t}D^-u(Y_s)dY_s + \int_{-\infty}^{\infty}\Lambda_t(a)\mu(da),
\end{align}
where $D^-u$ is the left-hand derivative of $u$, $\Lambda(a)$ is the local time of $Y$ in $a$, and $\mu$ is the second derivative measure, i.e. 
$\mu([a,b))=D^-u(b)-D^-u(a)$, $a<b$ real. 
Notice that the first derivative  of $u$ exists and is continuous outside $L-5R$, and the second derivative of $u$ exists (in
particular) outside the Lebesgue zero set $A= \{-bL+5R,0,L-5R\}$.
Hence we find by definition of the second derivative measure, and  with the help of equation (7.3) in \cite{kar-shr} that $P'_{0,\tilde V_{\ell,b,L}}$-a.s.,
\begin{equation}
\label{eq:sec-der}
\begin{aligned}
\int_{-\infty}^{\infty}\Lambda_t(a)\mu(da)=& \int_{-\infty}^{\infty}
\Lambda_t(a){\bf 1}_{A^c}(a) u''(a)\,da +\Lambda_t(L-5R)\,j(L-5R)\\
=& \tfrac{1}{2}\int_0^t  u''(Y_s){\bf 1}_{A^c}(Y_s) d\langle Y \rangle_s+
\Lambda_t(L-5R)\,j(L-5R)\,.
\end{aligned}
\end{equation}
Another application of  equation (7.3) p.218 in \cite{kar-shr} shows that 
\begin{equation}
\label{eq:N-0}
P'_{0,\tilde V_{\ell,b,L}}-\text{a.s.,}\quad  \int_0^t {\bf 1}_{A}(Y_s)d\langle Y \rangle_s
=2\int_{-\infty}^\infty {\bf 1}_{A}(a)\Lambda_t(a)\,da=0\,.
\end{equation}
Since $d\langle Y \rangle_s= 
 l \cdot a'_{\tilde V_{\ell,b,L}}(X_s)\,l \,ds 
\ge ds/\nu$,
see Lemma \ref{lemma:L'}, we deduce from 
(\ref{eq:N-0}) that, 
\begin{equation}
  \label{eq:N}
  P'_{0,\tilde V_{\ell,b,L}}-\text{a.s.,}\,\,   \int_0^t {\bf 1}_{A}(Y_s)\,ds=0\,,\text{ and hence, }
  P'_{0,\tilde V_{\ell,b,L}}-\text{a.s.,}\,\,   \int_0^t D^-u(Y_s){\bf 1}_{A}(Y_s)dY_s =0\,.
\end{equation}
Combining (\ref{eq:sec-der}) and (\ref{eq:N}), and by definition of the operator $\mathcal L'$, see below (\ref{eq:b'}), 
we can now rewrite (\ref{eq:ito}) as the
$P'_{0,\tilde V_{\ell,b,L}}$-a.s. equalities
\begin{equation*}
\begin{aligned}
u(Y_{t})=&1+\int_0^{t} u'(Y_s){\bf 1}_{A^c}(Y_s)dY_s 
+\tfrac{1}{2}\int_0^{t} u''(Y_s)\, {\bf 1}_{A^c}(Y_s)d\langle Y \rangle_s + \Lambda_t(L-5)j(L-5R)\\
=&1+\int_0^{t}\mathcal L'\tilde u(X_s){\bf 1}_{A^c}(Y_s)\,ds + \Lambda_t(L-5R)\,j(L-5R)
+M_t\,,
\end{aligned}
\end{equation*}
where $M_t$ is a continuous martingale.
In particular, 
$\tilde u (X_t)$ (=$u(Y_t)$) is a continuous semimartingale, and 
applying It\^o's rule to the product $e^{\lambda t} \cdot \tilde u(X_t)=\tilde v_\lambda(t,X_t)$, and using (\ref{eq:N}) once again, we obtain that,
$P'_{0,\tilde V_{\ell,b,L}}$-a.s.,
\begin{equation}
\label{eq:ito-2}
\begin{aligned}
&v_{\lambda}(t,X_{t })
=\,1+\int_0^{t }\lambda e^{\lambda s}
  \tilde u(X_s) \,ds + \int_0^{t } e^{\lambda s}\, d\tilde u(X_s)\\
&=\,1+ \int_0^{t } \Big(\tfrac{\partial}{\partial s}+
   \mathcal L'\Big)\tilde v_\lambda (s,X_s){\bf 1}_{A^c}(Y_s)\,ds
 +j(L-5R)\int_0^{t } e^{\lambda s}d\Lambda^{L-5R}_s +N_{t }\,,
\end{aligned}
\end{equation}
where $N_{t }$ is a continuous martingale.
Using $\tfrac{1}{\nu} \le l \cdot a'_{\tilde V_{\ell,b,L}}(x)\,l \le \nu$, we find through direct computation that for $ x \in V_{\ell,b,L}$, and a suitable $\psi (x) \geq 0$, using the notation $I_1=(-bL,-bL+5R)$, $I_2=(-bL+5R,L-5R)$, $I_3=(L-5R,L)$,
\begin{equation*}
 [(\tfrac{\partial}{\partial s}+ \mathcal L')\tilde v_\lambda](s,x)\leq \psi (x)e^{\lambda s} \cdot
\begin{cases}
\lambda(e^{20\nu \beta R}\alpha_3-1)-4\nu \beta( 2 \beta + b'_{\tilde V_{\ell,b,L}} (x) \cdot \ell)\,,&\text{ if } x \cdot \ell \in I_1\,,\\ 
\lambda+\alpha_2\tfrac{\epsilon}{\nu}(\frac{1}{2}\alpha_2\epsilon-b'_{\tilde V_{\ell,b,L}} (x) \cdot \ell)\,,&\text{ if }x \cdot \ell \in I_2\,,\\
\lambda(\alpha_5-1)-4\nu \beta( 2 \beta + b'_{\tilde V_{\ell,b,L}} (x) \cdot \ell)\,,&\text{ if }x \cdot \ell \in I_3\,.
\end{cases}
\end{equation*}
Hence, by (\ref{eq:b'-bounds}) and (\ref{eq:coefficients}), we can find 
$\lambda_0>0$ small such that for $x \in V_{\ell,b,L}$, $x \cdot \ell \notin A$, the right-hand side of the last expression is negative. 
Since $j(L-5R) \le 0$, see (\ref{eq:derivative}), we obtain from (\ref{eq:ito-2}) applied to the stopping time $t \wedge T_{V_{\ell,b,L}}$ that (\ref{eq:supermartingale}) holds.\\ 
We now derive the claim of the proposition from (\ref{eq:supermartingale}).
When $d \geq 2$, the probability to exit $V_{\ell,b,L}$ neither from the ``right'' nor from the ``left''
can be bounded as follows:
\begin{equation}
\label{eq:bound-0}
\begin{aligned}
  & P'_{0,\tilde V_{\ell,b,L}}[-bL<X_{T_{V_{\ell,b,L}}} \cdot \ell < L\,]\leq \\
 & P'_{0,\tilde V_{\ell,b,L}}[-bL<X_{T_{V_{\ell,b,L}}} \cdot \ell < L,\,T_{V_{\ell,b,L}} > 
\tfrac{2\alpha_2 \epsilon}{\lambda_0}L\,]
+ P'_{0,\tilde V_{\ell,b,L}}[\,\sup |X_t|\geq L^2:t \leq \tfrac{2\alpha_2 \epsilon}{\lambda_0}L\,]\,.
\end{aligned}
\end{equation}
By Chebychev's inequality and Fatou's lemma, we find that the first term on the right-hand 
side is smaller than
\begin{equation}
\label{eq:bound-1}
\begin{aligned}
&\frac{1}{v_{\lambda_0}(\tfrac{2\alpha_2 \epsilon}{\lambda_0}L,L)}E'_{0,\tilde V_{\ell,b,L}}
[v_{\lambda_0}(T_{V_{\ell,b,L}},X_{T_{V_{\ell,b,L}}}\cdot \ell)]\\
\le& c(\epsilon) e^{-\alpha_2 \epsilon L} \liminf_{t \to \infty}E'_{0,\tilde V_{\ell,b,L}}
[v_{\lambda_0}(t \wedge T_{V_{\ell,b,L}},X_{t \wedge T_{V_{\ell,b,L}}}\cdot \ell)]\\
\le& c(\epsilon) e^{-\alpha_2 \epsilon L} v_{\lambda_0}(0,0)= c(\epsilon) e^{-\alpha_2 \epsilon L},
\end{aligned}
\end{equation}
where, in the last inequality, we used (\ref{eq:supermartingale}). 
Applying Lemma \ref{lemma:bernstein} iii) to the second term in the right-hand 
side of (\ref{eq:bound-0}), we obtain together with (\ref{eq:bound-1}), that
\begin{equation}
\label{eq:bound-2}
\limsup_{L \to \infty}L^{-1} \log P'_{0,\tilde V_{\ell,b,L}}[-bL < X_{T_{V_{\ell,b,L}}}\cdot \ell<L]
<0\,.
\end{equation}
When $d \geq 1$, we bound the probability to exit $V_{\ell,b,L}$ from the left
by a similar argument as in (\ref{eq:bound-1}), and find that
\begin{equation}
  \label{eq:bound-3}
   P'_{0,\tilde V_{\ell,b,L}}[X_{T_{V_{\ell,b,L}}} \cdot \ell=-bL] \leq \frac{v_{\lambda_0}(0,0)}{v_{\lambda_0}(0,-bL)}
\leq e^{-c(\epsilon)L}\,.
\end{equation}
(\ref{eq:bound-3}), together with  (\ref{eq:bound-2}), when $d \geq 2$, show 
(\ref{eq:exit-bounded-set}), which implies condition $(T)|\ell$.
\end{proof}

Now the {\bf proof of Theorem \ref{thm:examples}} is carried out by checking
condition $(K)|\ell$. The proof remains the same as for the case $a=
Id$, 
which can be found in \cite{schmitz}, see (5.38) and the subsequent computations.\\
This finishes the proof of Theorem \ref{thm:examples}. 
$\hfill{}\Box$

\section{Condition $(T)$ for a diffusion with divergence-free drift}

In this section, we assume that $d \geq 2$. On $(\Omega,\mathcal A, \mathbb P)$ we consider a random variable $H(\cdot)$ with values in the space of skew-symmetric $d \times d$-matrices, and we define for $\om \in \Omega,x \in \R{d}$, $H(x,\om) = H(t_x \om)$, 
so that $H$ is stationary under $\mathbb P$.
We assume that $H$ obeys finite range dependence, i.e.
\begin{equation}
\label{eq:H-1}
(\ref{eq:R-separation}) \text{ holds with $H$ in place of $a$, $b$},
\end{equation}
and that $H$ is bounded and has bounded and Lipschitz continuous 
derivatives, i.e. 
there are constants $\eta,\,\tilde \beta,\,\tilde K$ such that
for all $1 \leq i,j,k \leq d$, $\om \in \Omega$,
\begin{align}
\bullet \,\,&|H(\cdot,\om)| \leq \eta\,\label{eq:H-2}\\
\bullet \,\, & \partial_{k} H_{ij}(\cdot,\om)\text{ exists, } |\partial_{k}H_{ij}(\cdot,\om)| \leq  \tilde \beta, \text{ and } \label{eq:H-3}\\
\bullet \,\,&|\partial_{k}H_{ij}(y,\om)-\partial_{k}H_{ij}(z,\om)| \leq  \tilde K |y-z|, \text{ for all }y,z \in \R{d}. 
\label{eq:H-4}
\end{align}
For a constant vector $v \neq 0$, we define diffusion and drift coefficients as
\begin{equation}
  \label{eq:b-div}
 a(\cdot,\om)= Id,\,\,\, b(\cdot,\om)= \nabla \cdot H(\cdot,\om)+v,\,\,\all\,,
\end{equation}
i.e. for $1 \leq j \leq d$, $b_j(\cdot,\om) = \sum_{i=1}^d \partial_{i} H_{ij}(\cdot,\om) 
+ v_j\,$. 
It follows from the skew-symmetry of $H$ and from (\ref{eq:b-div}) that for all $\om$, $\nabla \cdot b(\cdot, \om) =0$ in the weak sense.  Then the measures
$P_{x,\om}$ can be viewed as describing the motion of a diffusing particle in  an incompressible fluid 
and equipped 
with a random stationary velocity field obeying finite range dependence. It further follows from (\ref{eq:b-div}) that the operator $\mathcal L_\om$ defined in (\ref{eq:diff-operator})  can be written with a principal part in divergence form: 
\begin{equation}
  \label{eq:div-form}
  \mathcal L_\om = \nabla \cdot ((\tfrac{1}{2}\, Id + H(\cdot,\om))\nabla) + v \nabla\,.
\end{equation}
\begin{thm}
\label{thm:div}
  ($d \geq 2$)\\
Assume (\ref{eq:H-2}), (\ref{eq:H-3}) and (\ref{eq:b-div}). 
Then $(T)|\hat v$ holds, where $\hat v = \tfrac{v}{|v|}$.
Moreover we have a strong law of large numbers with an explicit velocity:
\begin{equation}
 \label{eq:velocity1}
 \text{for all }\all, \,\,\,P_{0,\om}-\text{a.s., }\,\, \lim_{t \to \infty}\frac{X_t}{t} \longrightarrow \mE[b(0,\om)]=v\,.
\end{equation}
If (\ref{eq:H-1})-(\ref{eq:b-div}) hold, then in addition a functional CLT holds:
\begin{equation}
\label{eq:clt1}
\begin{aligned}
&\text{ under $P_{0}$, $B_{\cdot}^{s}=\tfrac{X_{\cdot s}-\cdot s \mE[b(0,\om)]}{\sqrt{s}}$
 converges in law on $C(\R{}_{+},\R{d})$, as $s \to \infty$, to a}\\
&\text{Brownian motion $B_{\cdot}$ with non-degenerate covariance matrix}.
\end{aligned}
\end{equation}
\end{thm} 
\begin{proof}
We recall the convention on constants stated at the end of the Introduction. It follows from (\ref{eq:H-3})-(\ref{eq:b-div}) that our standing assumptions (\ref{eq:a}),\,(\ref{eq:Lipschitz}) hold with $\nu=1$, $\beta=d^2 \tilde \beta+d|v|$ and $K=d^2 \tilde K$.
The following upper bound on the heat kernel from Theorem 1.1 in Norris \cite{norris} will be instrumental in checking $(T)|\hat v$ and in proving (\ref{eq:velocity1}): 
there is a constant $c_2$ that depends only on $\eta$ and $d$, such that for all $\om \in \Omega$, $x,y \in \R{d}$ and $t>0$,
\begin{equation}
  \label{eq:heat_kernel}
  p_\om(t,x,y+vt) \leq c_2\, t^{-d/2}\exp\{-\tfrac{|y-x|^2}{c_2\,t}\}\,.
\end{equation}
We first prove (\ref{eq:velocity1}). 
Choose $\varepsilon>0$. For $n \ge 1$, we define $A_n=\{|\tfrac{X_n}{n}-v|>\varepsilon\}$ and $B_n=\{\tfrac{Z_1 \circ \theta_n}{n}>\varepsilon\}$
(recall $Z_1$ in Lemma \ref{lemma:bernstein}).
With the help of (\ref{eq:heat_kernel}), we find for $\all$, and for $n$ large,
\begin{align*}
P_{0,\om}[A_n]=\int_{B^c_{\varepsilon n}(0)}p_\om(n,0,y+nv)\,dy
\le  \int_{B^c_{\varepsilon n}(0)} c_2\,n^{-d/2} e^{-\tfrac{|y|^2}{c_2\,n}}\,dy
\le c_2^2\, n^{1-d/2}\, e^{-\tfrac{\varepsilon^2}{2c_2}\,n}\,.
\end{align*}
From the Markov property, and from Lemma \ref{lemma:bernstein} ii), we find that
\begin{align*}
\sup_\om P_{0,\om}[B_n]
\le \sup_{x,\om} P_{x,\om}[Z_1>\varepsilon n] \le \tilde c(\varepsilon)\, e^{-c(\varepsilon) n^2}\,.
\end{align*}
For $t>0$, we write 
\begin{align*}
\tfrac{X_t}{t}=\tfrac{\lfloor t \rfloor}{t}\,\left(\tfrac{X_{\lfloor t \rfloor}}{\lfloor t \rfloor}\,+\,\tfrac{X_t-X_{\lfloor t \rfloor}}{\lfloor t \rfloor}\right),
\end{align*}
and an application of Borel-Cantelli's lemma to the sets $A_n$ resp.~$B_n$ 
shows that for all $\all$, $P_{0,\om}$-a.s., $\lim_t \tfrac{X_t}{t}=v$.
Observe that, due to the stationarity of $b$ and $H$, for all $x \in \R{d}$,
\begin{align*}
  \mathbb E[b(0,\om)]=
\mathbb E[b(x,\om)]=\nabla \cdot \int_\Omega H(t_x\om)\,
  \mathbb P(d \om)\,+\,v =\nabla \cdot \int_\Omega H(\om)\,\mathbb P(d \om)\,+\,v  =v,
\end{align*}
and (\ref{eq:velocity1}) follows.
Once we have checked $(T)|\hat v$, then (\ref{eq:clt1}) immediately follows from (\ref{eq:clt}) and (\ref{eq:velocity1}). To conclude the proof of Theorem \ref{thm:div},
it thus remains to show that $(T)|\hat v$ holds.\\
For $x \in \R{d}$ and $\ell \in S^{d-1}$, define the projection on the orthogonal
complement of $\ell$, as well as  a bounded approximation $V_{\ell,b,L}$ of the slab $U_{\ell,b,L}$, 
\begin{equation}
\label{eq:V}
 \Pi_{\ell}(x) \df x-(x \cdot \ell)\,\ell\,,\quad   V_{\ell,b,L} \df \{x \in U_{\ell,b,L}: |\Pi_\ell(x)| < L^2\}\,\,.
\end{equation}
We choose unit vectors $\ell'$ such that $\hat v \cdot \ell'>\tfrac{2}{3}$, and  
we will show that for all such $\ell'$ and all $b>0$, 
\begin{equation}
\label{eq:exit-0}
  \limsup_{L \to \infty}L^{-1}\sup_{\all} \log\,P_{0,\om}[ X_{T_{V_{\ell',b,L}}} \cdot \ell' <L]<0\,,
\end{equation}
which clearly implies condition $(T)|\hat v$.
Define $\gamma = \tfrac{3}{|v|}$, and observe that for all $\om \in \Omega$,
\begin{equation}
\label{eq:exit-1}
 P_{0,\om}[X_{T_{V_{\ell',b,L}}} \cdot \ell' <L] \leq P_{0,\om}[X_{T_{V_{\ell',b,L}}} \cdot \ell' <L, T_{V_{\ell',b,L}} \leq \gamma L]+ P_{0,\om}[ T_{V_{\ell',b,L}} > \gamma L]\,.
\end{equation}
We first estimate the second term on the right-hand side. 
Notice that for all $y \in V_{\ell',b,L}$, $|y-\gamma Lv| \geq ((L\ell'-\gamma Lv)\cdot \ell')_-=L(3\hat v \cdot \ell'-1) > L$.
An application of (\ref{eq:heat_kernel}) yields that 
\begin{equation}
\label{eq:estimate-1}
\begin{aligned}
 &\sup_\om P_{0,\om}[T_{V_{\ell',b,L}} > \gamma L] \leq  \sup_\om P_{0,\om}[X_{\gamma L} \in V_{\ell',b,L}] \\
\leq& c_2 \gamma^{-d/2} L^{-d/2} \int_{V_{\ell',b,L}} \exp\{-\tfrac{|y-v\gamma L|^2}{c_2\gamma L}\}\,dy
 \leq c(b,\gamma,\eta)L^{3d/2-1} \exp\{-\tfrac{L}{c_2 \gamma}\}\,.
\end{aligned}
\end{equation}
In view of (\ref{eq:exit-0}) and (\ref{eq:exit-1}), it remains to show that for all $b>0$, 
and $\ell'$ as above (\ref{eq:exit-0}),
\begin{equation}
\label{eq:estimate3}
\limsup_{L \to \infty}L^{-1}\sup_{\all} \log\,P_{0,\om}[ X_{T_{V_{\ell',b,L}}} \cdot \ell' <L, T_{V_{\ell',b,L}} \leq \gamma L]<0.
\end{equation}
Introduce the subsets of $V_{\ell',b,L}$,
\begin{align*}
V_{\ell',b,L}^- \df \{x \in V_{\ell',b,L}:x \cdot \ell' <-\tfrac{bL}{2}\}\,,\,\,
  V_{\ell',b,L}^0 \df \{x \in V_{\ell',b,L}: |\Pi_{\ell'}(x)|>\tfrac{L^2}{2} \}\,,
\end{align*}
as well as the union of these two sets, $\tilde V_{\ell',b,L} \df V_{\ell',b,L}^- \cup V_{\ell',b,L}^0$, and write
\begin{align}
&P_{0,\om}[X_{T_{V_{\ell',b,L}}} \cdot \ell' <L, T_{V_{\ell',b,L}} \leq \gamma L]\,\, \leq \,\,
P_{0,\om}[T_{V_{\ell',b,L}}\leq 1]\,\,+ \label{eq:estimate-2}\\
+&\sum_{n=1}^{\lfloor \gamma L \rfloor}
P_{0,\om}[T_{V_{\ell',b,L}} \in (n,n+1],\,X_{T_{V_{\ell',b,L}}} \cdot \ell' <L,\,
X_n \notin \tilde V_{\ell',b,L}]+
\sum_{n=1}^{\lfloor \gamma L \rfloor}P_{0,\om}[X_n \in \tilde V_{\ell',b,L}]\,.\nonumber
\end{align}
Lemma \ref{lemma:bernstein} ii) shows that 
$\sup_\om P_{0,\om}[T_{V_{\ell',b,L}}\leq 1] \le \tilde c(b) \exp\{-c(b)L^2\}$.
With the help of the Markov property and Lemma \ref{lemma:bernstein} ii) 
we find that for large $L$, the middle term on the right-hand side of (\ref{eq:estimate-2}) is smaller than
\begin{equation}
\begin{aligned}
\sum_{n=1}^{\lfloor \gamma L \rfloor}\sup_\om P_{0,\om}[\sup_{0 \leq t \leq 1}
|X_t - X_0|\circ \theta_n>\tfrac{bL}{2}]
\leq& \,\gamma L \sup_{x,\om} P_{x,\om}[\sup_{0 \leq t \leq 1}|X_t - X_0|>\tfrac{bL}{2}]\\
\leq&\, \tilde c(b,\gamma) L \,\exp\{-c(b) L^2\}.
\end{aligned}
\end{equation}
We now bound the third term on the right-hand side of (\ref{eq:estimate-2}):
\begin{equation}
\label{eq:estimate-4}
 \sum_{n=1}^{\lfloor \gamma L \rfloor}\sup_\om P_{0,\om}[X_n \in \tilde V_{\ell',b,L}] 
 \leq \sum_{n=1}^{\lfloor \gamma L \rfloor} 
 \sup_\om P_{0,\om}[X_n \in  V_{\ell',b,L}^-]+
 \gamma L\, \sup_\om P_{0,\om}[\sup_{t \leq \gamma L}|X_t|>L^2/2]\,.
\end{equation}
Lemma \ref{lemma:bernstein} iii) shows that
 the second term on the right-hand side is smaller than
$\tilde c(\gamma) L \exp\{-c(\gamma)L^3\}$.  
Using (\ref{eq:heat_kernel}) we can bound the first term 
on the right-hand side of (\ref{eq:estimate-4}) with
\begin{equation}
  \label{eq:estimate-6}
  \sum_{n=1}^{\lfloor \gamma L \rfloor}\int_{V_{\ell',b,L}^-}c_2\,n^{-d/2}
\exp\{-\tfrac{|y-nv|^2}{c_2 n}\}\textrm{d}y
\leq c(b,\gamma,\eta)\, L^{2d}\,\exp\{-\tfrac{b^2}{4c_2\gamma}\,L\}\,,
\end{equation}
where, in the last inequality, we simply used that $1 \leq n \leq \gamma L$ and, and that for all $n \geq 1$, 
and $y \in V_{\ell',b,L}^-$, $|y-nv|\geq ((y-nv)\cdot \ell')_- \geq bL/2$ since
$\ell' \cdot \hat v >0$. 
We now obtain (\ref{eq:estimate3}) by
collecting the results between (\ref{eq:estimate-2}) and (\ref{eq:estimate-6}).This finishes the proof of the theorem.
\end{proof}

\section{Condition $(T)$ for an anisotropic gradient-type diffusion}

Let $\varphi$ be a real-valued random variable on $(\Omega,\mathcal A, \mathbb P)$, and define for $\om \in \Omega,x \in \R{d}$, $\varphi(x,\om) \df \varphi(t_x \om)$,
so that $\varphi$ is stationary.
We assume that $\varphi$ obeys finite range dependence, i.e.
\begin{equation}
\label{eq:V-1}
(\ref{eq:R-separation}) \text{ holds with $\varphi$ in place of $a$, $b$},
\end{equation}
and that $\varphi$ is bounded and has bounded and Lipschitz continuous derivatives, i.e. 
there are constants $\eta,\,\beta,\,K$ such that
for all $1 \leq i,j,k \leq d$, $\om \in \Omega$,
\begin{align}
\bullet \,\,&|\varphi(\cdot,\om)| \leq \eta\,\label{eq:V-2}\\
\bullet \,\, &|\nabla \varphi(\cdot,\om)|\le \tilde \beta, \text{ and }
|\nabla \varphi(y,\om)-\nabla \varphi(z,\om)|\le K|y-z|,\,\,y,\,z \in \R{d}.
\label{eq:V-3}
\end{align}
Fix $\lambda>0$  and $\ell \in S^{d-1}$, and define the non-stationary function $\psi(x, \om)= \varphi(x, \om)+\lambda \ell \cdot x$, 
$\all$, $x \in \R{d}$. 
We now  define diffusion and drift coefficients
\begin{equation}
  \label{eq:b-grad}
 a(\cdot,\om) = Id,\quad b(\cdot,\om)= \nabla \psi(\cdot,\om)=\nabla \varphi(\cdot,\om)+\lambda \ell\,,
  \,\, \text{ all }\all\,.
\end{equation}
It follows from (\ref{eq:V-2}) that
\begin{equation}
  \label{eq:phi}
  e^{-2\eta}e^{2 \lambda \ell \cdot x} \leq e^{2\psi(x,\om)} \leq e^{2\eta}e^{2 \lambda \ell \cdot x}\,,\quad x \in \R{d},\,\all.
\end{equation}
Under the above assumptions, Shen has shown that a ballistic
law of large numbers and an invariance principle holds under the annealed measure $P_0$, see Section 4 of \cite{shen}. Our aim is to verify that condition $(T)|\ell$ holds.\\

\begin{thm}
\label{thm:grad}
$(d \geq 1)$ Condition $(T)|\ell$ holds under the assumptions (\ref{eq:V-2})-(\ref{eq:b-grad}). 
\end{thm}
\begin{proof}
We recall the convention on constants stated at the end of the Introduction. 
It follows from (\ref{eq:V-3}),(\ref{eq:b-grad}) that our standing assumptions (\ref{eq:a}),\,(\ref{eq:Lipschitz}) hold with $\nu=1$, $\beta= \tilde \beta+\lambda$ and $K$.
We first introduce some notation. Recall the projection $\Pi$ in (\ref{eq:V}), and define for $\ell' \in S^{d-1}$ and $\delta >0$, $V_{\ell',\delta}=\{x \in U_{\ell',b,L}:|\Pi_{\ell'}(x)| < L/\sqrt{\delta}\}$.
In fact, we will show that there exists $0<\delta<1$ such that for all $\ell ' \in S^{d-1}$ with 
$\ell' \cdot \ell  > \sqrt{1-\delta^2}$,
\begin{equation}
  \label{eq:exit}
   \limsup_{L \to \infty}L^{-1} \sup_\om \log P_{0,\om}[ X_{T_{V_{\ell',\delta}}} \cdot \ell' <L]<0\,,
\end{equation}
which easily implies condition $(T)|\ell$.
The strategy of proof of (\ref{eq:exit}) relies on the methods introduced in Section 4.1 in Shen \cite{shen}. We recall some important facts from there.\\
We denote with $P_\om^t f(x)=E_{x,\om}[f(X_t)]$, $f$ bounded measurable, the quenched 
semi-group generated by the operator $\mathcal L_\om$ defined in (\ref{eq:diff-operator}), and with $m_\om(dx)=\exp(2\psi(x,\om))\,dx$ the reversible measure to $P^t_{\om}$. 
The following key estimate is contained in Proposition 4.1 of \cite{shen}: there is a positive constant $c_3(\eta,\lambda)$ such that
\begin{equation}
  \label{eq:spectral-gap}
  \sup_{\om, U}\|P^t_{\om,U}\|_{m_\om} \leq \exp\{-c_3 t\}\,, \, t>0\,,
\end{equation}
 where $(P^t_{\om,U}f)(x) \df E_{x,\om}[f(X_t), T_U >t]$, $t>0$, $f \in  L^2(m_\om)$, 
$\|\cdot\|_{m_\om}$ denotes the operator norm in 
$L^2(m_\om)$ and $U$ varies over the collection of non-empty open subsets of $\R{d}$. \\
The claim (\ref{eq:exit}), and hence Theorem \ref{thm:grad}, follow from the next two propositions:
\begin{prop}
  \label{prop:1}
Let $0<\delta <1$. There is a positive constant $c_4(\delta,b,\eta)$ such that for $\ell' \in S^{d-1}$ with $\ell' \cdot \ell > \sqrt{1-\delta^2}$
and $L>0$,
\begin{equation}
  \label{eq:exit-time}
  \sup_\om P_{0,\om}[\,T_{V_{\ell',\delta}} \geq \tfrac{3\lambda}{c_3} L] 
  \leq c_4 e^{-\tfrac{\lambda}{2} L}\,.
\end{equation}
\end{prop}
\begin{prop}
  \label{prop:2}
There exist $0<\delta<1$ and positive constants $c_5(\delta,b,\eta,\lambda),\,c_6(b,\eta,\lambda)$ such that for all $\ell ' \in S^{d-1}$ with $\ell' \cdot \ell  > \sqrt{1-\delta^2}$ and $L>0$,
\begin{equation}
  \label{eq:exit1}
\sup_\om P_{0,\om}[T_{V_{\ell',\delta}} < \tfrac{3\lambda}{c_3} L,\,\,
 X_{T_{V_{\ell',\delta}}} \cdot \ell' <L] \leq c_5 e^{-c_6 L}\,.   
\end{equation}
\end{prop}
The proofs are close to the proofs of Propositions 4.2 and 4.3 in Shen 
\cite{shen}.  We indicate the main steps of the proofs in the Appendix.
\end{proof}

\section{Appendix}

\subsection{Bernstein's Inequality}
\label{sec:appendix}
The following Lemma follows in essence from Bernstein's inequality (see 
\cite {rev-yor} page 153-154). A proof of an inequality similar to those 
below can be found in the appendix of \cite{schmitz}. 
Recall the definition of the class of operators $\cN(\nu,\beta)$ at the beginning 
of Section \ref{sec:2}.
\begin{lemma}
\label{lemma:bernstein}
Let $\cL \in \cN(\nu,\beta)$, denote with
$Q_x$ an arbitrary solution to the martingale problem for $\cL$ started at 
$x \in \R{d}$, and set 
$Z_t = \sup_{s \leq t}|X_s-X_0|$.
For every $\gamma >0$,
there are positive constants $c,\,\tilde c$, depending only on $\gamma,\,\nu,\,\beta,\,d$ such that for $L>0$, 
\begin{equation}
\label{eq:bern}
\begin{aligned}
i)\,\,&\sup_x Q_x \big[Z_{\gamma L}\geq 2\gamma \beta L\big]\leq \tilde c\,e^{-c\,L},\quad 
ii)\,\, \sup_x Q_x \big[Z_1 \geq \gamma L\big]\leq \tilde c\,e^{-c\,L^2},\\
iii)\,\,&\sup_x Q_x\big[Z_{\gamma L}\geq L^2\big] \le \tilde c\,e^{-c\,L^3}\,.
\end{aligned}
\end{equation}
\end{lemma}

\subsection{Bounds on the transition density and on the Green function}
Recall the convention on the constants stated at the end of the Introduction. 
\begin{lemma}
  \label{lemma:PDE}
Let $\mathcal L_\om$ be as in (\ref{eq:diff-operator}), and let assumptions
(\ref{eq:a})
and (\ref{eq:Lipschitz}) be in force. Then the linear parabolic equation of second order $\frac{\partial u}{\partial t}=\mathcal L_\om u$ has a   unique fundamental solution $p_\om(t,x,y)$, and there are positive constants $c,\tilde c$, such that for $\all$, 
and $0<t\le 1$,
\begin{equation}
    \label{eq:PDE-upper}
    |p_\om(t,x,y)| < \tfrac{\tilde c}{t^{d/2}} \exp\big\{-\tfrac{c |x-y|^2}{t}\big\}\,.
  \end{equation}
\end{lemma}
For the proof of (\ref{eq:PDE-upper}) we refer the reader to \cite{illin}. The statement (4.16) therein corresponds to (\ref{eq:PDE-upper}). 
Recall the Green function $g_{\om,U}$ in (\ref{eq:green-function}). In Corollary 6.3 in \cite{schmitz}, we obtained 
the following
\begin{cor}
 \label{cor:green}
Assume (\ref{eq:a}) and (\ref{eq:Lipschitz}),
 and let $U$ be a bounded $C^\infty$-domain. 
There are positive constants $\tilde c, c(U)$ such that for $x \in U$, and all
$\om \in \Omega$,
\begin{equation}
  \label{eq:green-upper}
  g_{\om,U}(x) \leq 
\begin{cases}
\tilde c\, |x|^{2-d}+ c, & \text{if $d \geq 3$ and $x \neq 0$},\\
\tilde c \log \frac{diam\,U}{|x|}+ c,& \text{if $d=2$ and $x \neq 0$},\\
c, &\text{if $d=1$}\,.
\end{cases}
\end{equation}
\end{cor}
\begin{proof}
We repeat here the computations from \cite{schmitz} for the convenience of the reader.
From the definition of $g_{\om,U}$ and $p_{\om,U}$, see (\ref{eq:green-function}) 
and (\ref{eq:duhamel}), we obtain
\begin{align}
\label{eq:sup-green}
  g_{\om,U}(x)=\int_{0}^{\infty}p_{\om,U}(t, 0,x)\,dt
   \leq \int_{0}^{1}p_{\om}(t, 0  ,x)\,dt + \sum_{k=2}^{\infty}\int_\frac{k}{2}^{\frac{k+1}{2}}p_{\om,U}(t,0,x)\,dt\,.
\end{align}
With the help of (\ref{eq:PDE-upper}), we find positive constants $\tilde c,\,c$
 such that $\int_{0}^{1}p_{\om}(t, 0  ,x)\,dt$ is smaller than the 
right-hand side of (\ref{eq:green-upper}), and hence, it suffices to show that for 
some constant $c(U)$, 
\begin{equation}
\label{eq:sum}
\sup_\om \sum_{k=2}^{\infty}\int_\frac{k}{2}^{\frac{k+1}{2}}p_{\om,U}(t,0,x)\,dt \le c <\infty.
\end{equation}
We obtain by a repeated use of the Chapman-Kolmogorov equation and by (\ref{eq:PDE-upper}), that for $k \geq 2$,
\begin{multline}
\label{eq:chap-kolm}
  \int_\frac{k}{2}^{\frac{k+1}{2}}p_{\om,U}(t,0,x)\,dt
  \leq  \int_U p_{\om,U}(\tfrac{1}{2},0,v)\,dv\,\, ~\sup_{v \in U}\int_\frac{  k-1}{2}^{\frac{k}{2}}p_{\om,U}(t,v,x)\,dt\\
  \stackrel{induction}{\leq}  \left(\sup_{v \in U}P_{v,\om}[T_{U}>\tfrac{1}{2}]  \right)^{k-1}
  \sup_{v \in U} \int_{\frac{1}{2}}^{1}p_{\om,U}(t,v,z)\,dt
  \leq c \left(\sup_{v \in U}P_{v,\om}[T_{U}>\tfrac{1}{2}]\right)^{k-1}.
\end{multline}
The Support Theorem of Stroock-Varadhan, see \cite{bass2} p.25, shows that\\ 
$\sup_{\om,\,v \in U}P_{v,\om}[T_{U}>\tfrac{1}{2}] 
\le 1-c(U)$ for some positive constant $c(U)$, and (\ref{eq:sum}) follows.
\end{proof}

\subsection{The proof of propositions \ref{prop:1} and \ref{prop:2}}
\begin{proof}[Proof of Proposition \ref{prop:1}]
For sake of notations, we write $V$ and $B$ for $V_{\ell',\delta}$ resp. $B_{\frac{L(1 \wedge b)}{2}}(0)$.
The Markov property and Cauchy-Schwarz's inequality yield for $t>0$
\begin{align*}
 P_{0,\om}[T_V >t]
\leq& E_{0,\om}[P_{X_1,\om}[T_V >t-1],X_1 \in B]+P_{0,\om}[X_1 \notin B]\\
\leq& \|1_B (\cdot) p_\om(1,0,\cdot)e^{-2\psi(\cdot,\om)}\|_{L^2(m_\om)}\,\|P^{t-1}_{\om,V}\|_{m_\om}\,\|1_V\|_{L^2(m_\om)}+P_{0,\om}[X_1 \notin B]\,,
\end{align*}
where $p_\om(s,x,y)$ denotes the transition density function under the quenched measure $P_{x,\om}$.
Lemma \ref{lemma:bernstein} ii) shows that $P_{0,\om}[X_1 \notin B] \le c(b) \exp (-c(b) L^2)$. Further, using (\ref{eq:phi})
and Lemma \ref{lemma:PDE}, we obtain 
\begin{equation*}
\|1_B (\cdot) p_\om(1,0,\cdot)e^{-2\psi(\cdot,\om)}\|_{L^2(m_\om)}^2 
\leq c \,e^{2\eta}\,\int  1_B (y)e^{-2\lambda \ell \cdot y} dy
\leq c(\eta, b)\,L^d e^{\lambda (1\wedge b)L}\,.
\end{equation*}
Choose a unit vector $\tilde \ell$, orthogonal to $\ell'$, and such that $\tilde \ell$ lies 
in the span of $\ell$ and $\ell'$.
Then $\ell \cdot \ell' > \sqrt{1-\delta^2}$ implies that $\ell \cdot \tilde
\ell <\delta$, and for $y \in V$, we find that 
$y \cdot \ell \leq (1+\sqrt{\delta})L <2L$. As a result, and using again (\ref{eq:phi})
and Lemma \ref{lemma:PDE},
\begin{equation}
\label{eq:5.3}
  \|1_V\|_{L^2(m_\om)}^2 \leq e^{2\eta}\, \int  1_V (y)e^{2\lambda \ell \cdot y} dy
\leq \,c(\delta,\eta) L^d \exp\{4 \lambda L\}.
\end{equation}
Jointly with (\ref{eq:spectral-gap}), we find that $ P_{0,\om}[T_V >t] \leq c(\delta,\eta,b)L^d \exp(\tfrac{5}{2} \lambda L)\exp(-c_3 t)$,
and the claim follows from the choice of $t$.
\end{proof}
\begin{proof}[Proof of Proposition \ref{prop:2}]
We write $V$ for $V_{\ell',\delta}$, $\ell' \in S^{d-1}$, and  $B$ for $B_{\frac{(1 \wedge b)L}{4}}(0)$,
and we recall the projection $\Pi$ in (\ref{eq:V}).
We define 
\begin{equation*}
  V^0 \df \{x \in V: |\Pi_{\ell'}(x)|>\tfrac{1}{2\sqrt{\delta}}L\}, \quad 
  V^- \df \{x \in V:x\cdot \ell' < -\tfrac{3}{4}(1\wedge b)L\}\,.
\end{equation*}
We proceed in a similar fashion as in the proof of Proposition \ref{thm:div}, see (\ref{eq:estimate-2}) and the following lines, and find with $\mu_0=3\lambda/c_3$ that 
\begin{equation*}
P_{0,\om}[T_V < \mu_0 L,\,\,
 X_{T_V} \cdot \ell' <L] \leq  
\sum_{n=1}^{\lfloor \mu_0 L \rfloor} P_{0,\om}[X_n \in V^- \cup V^0]+c(\delta,b) \exp\{-c(\delta,b)L\}\,.
\end{equation*}
The  first term in the right-hand side can be bounded by
\begin{equation}
\label{eq:1}
 \sum_{n=1}^{\lfloor \mu_0 L \rfloor} P_{0,\om}[X_n \in V^0] + \sum_{n=1}^{\lfloor \mu_0 L \rfloor}P_{0,\om}[X_n \in V^-, X_1 \in B] +\mu_0L\, P_{0,\om}[X_1 \notin B]\,.
\end{equation}
Lemma \ref{lemma:bernstein} ii) shows that 
the righmost term in (\ref{eq:1}) is smaller than $c(b,\eta,\lambda)L e^{-c(b)L}$, and, when in 
addition $\delta^{-1/2}\ge 4\mu_0\beta$, Lemma \ref{lemma:bernstein} i) 
shows that the leftmost sum in (\ref{eq:1}) is smaller than $c(\eta,\lambda)L\exp\{-c(\eta,\lambda)L\}$.
It thus suffices to bound the middle term in (\ref{eq:1}).
Introduce the measure $m(dx)=\exp\{2\lambda \ell \cdot x\}dx$. It follows from the Markov property and Lemma \ref{lemma:PDE} that the second sum in (\ref{eq:1}) is smaller than
\begin{equation}
\label{eq:2}
  \sum_{n=1}^{\lfloor \mu_0 L \rfloor} \int_B p_\om(1,0,y)(P^{n-1}_\om 1_{V^-})(y)dy
\leq \sum_{n=1}^{\lfloor \mu_0 L \rfloor} c \exp\{\tfrac{1}{2}\lambda (1 \wedge b)L\} \langle P^{n-1}_\om 1_{V^-},1_B\rangle_m \,,
\end{equation}
where, for two measurable functions $f,g$, $\langle f,g \rangle_m$ denotes the integral of the product $f \cdot 
g$ w.r.t. the measure $m$.
It follows from Theorem 1.8 in \cite{sturm} (for details see the proof of Proposition 4.3 in \cite{shen})
that there is a constant $c_7(\eta)$ such that for all $\all$, and any open sets $A,D$ in $\R{d}$,
\begin{equation}
\label{eq:sturm}
\langle P^{n-1}_\om 1_A,1_D\rangle_m \leq \sqrt{m(A)}\sqrt{m(D)}
\exp\{-\tfrac{d(A,D)^2}{4c_7(n-1)}\} \,.
\end{equation}
We find by definition of $B,\,V$ and, by a similar argument as given 
above (\ref{eq:5.3}), that
\begin{equation*}
  m(B) \leq c(b) L^d \exp\{\tfrac{1}{2}\lambda(1 \wedge b) L\},\,\,\,
  m(V^-) \leq c(\delta,b) L^d \exp\{-\tfrac{3}{2}\lambda(1 \wedge b)\sqrt{1-\delta^2}\,L+2\lambda\sqrt{\delta} L\}. 
\end{equation*}
Since $d(B,V^-)=\tfrac{1}{2}(1 \wedge b)L$, 
(\ref{eq:sturm}) applied to the open sets $B,V$, yields that for large $L$, the right-hand side of (\ref{eq:2}) is 
smaller than
\begin{align*}
&\sum_{n=1}^{\lfloor \mu_0 L \rfloor}c \exp\{\tfrac{1}{2}\lambda (1 \wedge b)L\} 
\sqrt{m(B)}\sqrt{m(V^-)}\,\exp\{-\tfrac{d(B,V_-)^2}{4c_7\mu_0L}\} \\
\leq&\,
c(\delta,b) \mu_0 L^{d+1}\exp\{\sqrt{\delta}\lambda L\}\,
\exp\{\tfrac{3}{4}\lambda(1 \wedge b)(1-\sqrt{1-\delta^2})L\}\,
\exp\{-\tfrac{(1 \wedge b)^2}{c(\lambda,\eta)}L\} \\
\leq &\,c(\delta,b)\, e^{- c(b,\eta,\lambda)L}\,,
\end{align*}
provided $\delta$ is chosen small enough. This finishes the proof.
\end{proof}


\begin{thebibliography}{30}\small

 
\bibitem{bass2}Bass, R.: ``Diffusions and Elliptic Operators'', Springer Verlag,  1998.

\bibitem{boboc-mustata}Boboc, N., Mustata, P.: ``Espaces harmoniques associ\'es aux op\'erateurs differentiels lin\'eaires du second ordre de type elliptique'', {\it Lecture Notes in Math.}, volume 68, Springer Verlag, 1968.




\bibitem{fan-kom}Fannjiang, A., Komorowski, T .: ``An invariance principle for diffusion in turbulence '', {\it Ann. Probab.}, {\bf 27}(2), page 751--781, 1999.

\bibitem{fan-papa}Fannjiang, A., Papanicolaou, G.: ``Diffusion in turbulence'', {\it Probab. Theory relat. Fields}, {\bf 105}, page 279--334, 1996.

\bibitem{fried}Friedman, A.: ``Stochastic Differential Equations and 
Applications, Vol.1, Academic Press, San Diego, 1975.

\bibitem{gil-tru}Gilbarg, D., Trudinger, N.S.: ``Elliptic Partial Differential Equations of the Second Order'', Springer Verlag, 1998.

\bibitem{illin}Il'in, A.M., Kalashnikov, A.S., Oleinik, O.A.: ``Linear equations of the second order of parabolic type'', {\it Russian Math. Surveys}, {\bf 17}(3), page 1--143, 1962.


\bibitem{kar-shr}Karatzas, I., Shreve, S.: ``Brownian Motion and Stochastic Calculus'', Second Edition, Springer Verlag, 1991.


\bibitem{kom-olla-01}Komorowski, T., Olla, S.: ``On homogenization of time-dependent random flows'',  {\it Probab. Theory relat. Fields}, {\bf 121}(1), page 98--116, 2001.


\bibitem{krylov92}Krylov, N.V.: ``On one-point weak uniqueness for elliptic
equations'',  {\it Comm. Partial Differential Equations}, {\bf 17}(11-12),
page 1759--1784, 1992.


\bibitem{krylov04}Krylov, N.V.: ``On weak uniqueness for some diffusions with discontinuous coefficients'',
{\it Stochastic Proc. Appl.}, {\bf 113}(1), page 37--64, 2004.


\bibitem{landim-olla-yau}Landim, C., Olla, S., Yau, H.T.: ``Convection-diffusion equation with space-time ergodic random flow'', {\it Probab. Theory relat. Fields}, {\bf 112}, page 203--220, 1998.


\bibitem{nad96}Nadirashvili, N.: ``Nonuniqueness in the Martingale Problem and the Dirichlet Problem for Uniformly Elliptic Operators'', {\it Ann. Scuola Norm. Sup. Pisa Cl. Sci.}, {\bf 24}(4), page 537--550, 1997.

\bibitem{norris}Norris, J.R.: ``Long-time Behavior of Heat Flows'', {\it Archive for rational Mechanics and Analysis}, {\bf 140}(2), page 161-195, 1997.

\bibitem{oel}Oelschl\"ager, K.: ``Homogenization of a diffusion process in a divergence-free random field'', {\it Ann. Probab.}, {\bf 16}(3), page 1084--1126, 1988.

\bibitem{olla94}Olla, S.: ``Homogenization of diffusion processes in random fields'', Ecole Doctorale, Ecole Polytechnique, Palaiseau, 1994.

\bibitem{olla01}Olla, S.: ``Central limit theorems for tagged particles and for diffusions in random environment''. In: ``Milieux Al\'eatoires'', Panoramas et Synth\`eses, Num\'ero 12, Soci\'et\'e Math\'ematique de France, 2001.

\bibitem{osada}Osada, H.: ``Homogenization of diffusion processes with random stationary coefficients'', {\it Lecture Notes in Math.}, volume 1021, page 507--517, Springer, Berlin, 1983.

\bibitem{osada1}Osada, H.: ``Diffusion processes with generators
of generalized divergence form'',  J. Math. Kyoto Univ., {\bf 27}(4),
page 597--619, 1987.


\bibitem{rev-yor}Revuz, D., Yor, M.: ``Continuous Martingales and Brownian Motion'', 3rd Edition, Springer Verlag, Berlin, 1999.



\bibitem{schmitz}Schmitz, T.: ``Diffusions in random environment and ballistic behavior'', accepted 
for publication in {\it Ann. I. H. Poincar\'e PR}, in press,
available at http://dx.doi.org/, with doi:10.1016/j.anihpb.2005.08.003.

\bibitem{shen}Shen, L.: ``On ballistic diffusions in random environment'', {\it Ann. I. H. Poincar\'e}, PR 39(5), page 839--876, 2003.

\bibitem{str-var}Stroock, D.W., Varadhan, S.R.S.: ``Multidimensional Diffusion Processes'', Springer Verlag, New York, 1979.

\bibitem{str}Stroock, D.W.: ``Diffusion semigroups corresponding to uniformly elliptic 
divergence form operators'', in: Lecture Notes in Math., Vol. 1321,
Springer Verlag, Berlin, page 316--347, 1988.
 
\bibitem{sturm}Sturm, K.T.: ``Analysis on local Dirichlet Spaces-II. Upper Gaussian estimates for the fundamental solutions of parabolic equations'',  {\it Osaka J. Math.}, {\bf 32},  page 275 --312, 1995.

\bibitem{szn00}Sznitman, A.-S.: ``Slowdown estimates and central limit theorem for random walks in random environment'', {\it J. Eur. Math. Soc.}, {\bf 2},  page 93--143, 2000.

\bibitem{szn01}Sznitman, A.-S.: ``On a class of transient random walks in random environment'',  {\it Ann. Probab.},  {\bf 29}(2), page 723--764, 2001.

\bibitem{szn02}Sznitman, A.-S.: ``An effective criterion for ballistic behavior of random walks in random environment'', {\it Probab. Theory relat. Fields}, {\bf 122}(4), page 509--544, 2002.

\bibitem{szn03}Sznitman, A.-S.: ``On new examples of ballistic random walks in random environment'', {\it Ann. Probab.}, {\bf 31}(1), page 285--322, 2003.

\bibitem{szn04}Sznitman, A.-S.: ``Topics in random walks in random environment'', {\it ICTP Lecture Notes Series}, Volume XVII: School and Conference on Probability Theory, May 2004, 
available at www.ictp.trieste.it/$\sim$pub\_off/lectures.

\bibitem{szn-zer}Sznitman, A.-S., Zerner, M.P.W.: ``A law of large numbers for random walks in random environment'', {\it Ann. Probab.}, {\bf 27}(4), page 1851--1869, 1999.

\bibitem{szn-zeit}Sznitman, A.-S., Zeitouni, O.: ``An Invariance Principle for Isotropic Diffusions in Random Environment'', 
 {\it Invent. Math.}, {\bf 164}(3), page 455--567, 2006.
  
\bibitem{zeit}Zeitouni, O.: ``Random Walks in Random Environment'', {\it Lecture Notes in Mathematics}, volume 1837, page 190--312, Springer, 2004.

\end{thebibliography}
\end{document}